\documentclass[11pt,a4paper]{article}

\usepackage{header}
\usepackage{abbrev}

\begin{document}

\title{Perfectly matched layers for the Boltzmann equation: stability and sensitivity analysis}

\author{Marco Sutti\thanks{National Center for Theoretical Sciences, Mathematics Division, Room 406, 4F of Chee-Chun Leung Cosmology Hall, National Taiwan University, No. 1, Sec. 4, Roosevelt Road, Taipei 10617, Taiwan (R.O.C.) (\email{msutti@ncts.ntu.edu.tw}).}\hspace{2mm}\orcidlink{0000-0002-8410-1372}
~and Jan S. Hesthaven\thanks{Computational Mathematics and Simulation Science, EPFL SB, MATH MCSS, MA C2 652 (Bâtiment MA), Station 8, CH-1015 Lausanne, Switzerland (\email{Jan.Hesthaven@epfl.ch}).}}

\date{\today}


\maketitle

\begin{abstract}
We study the stability and sensitivity of an absorbing layer for the Boltzmann equation by examining the Bhatnagar--Gross--Krook (BGK) approximation and using the perfectly matched layer (PML) technique. To ensure stability, we discard some parameters in the model and calculate the total sensitivity indices of the remaining parameters using the ANOVA expansion of multivariate functions. We conduct extensive numerical experiments on two test cases to study stability and compute the total sensitivity indices, which allow us to identify the essential parameters of the model.

\bigskip
\textbf{Key words.} BGK model, perfectly matched layer, differential operators, stability analysis, ANOVA expansion, total sensitivity index

\medskip
\textbf{AMS subject classifications.} 76P05, 65N06, 35L45, 35B35, 35L05, 35Q30 

\end{abstract}

\section{Introduction}

Due to physical constraints on computational resources, numerical simulations of unbounded physical problems are virtually impossible to carry out without the truncation of the simulation domain. When the solution is not periodic, one must truncate the domain by introducing a boundary or absorbing layer.

In this work, we study and enhance an effective absorbing layer for the Boltzmann equation. In particular, we consider the Bhatnagar--Gross--Krook (BGK) approximation of the Boltzmann equation~\protect{\cite{BGK:1954}} in the approximate form proposed by~\protect{\cite{Tolke:2000}}, and use the perfectly matched layers (PML) technique.

Bérenger~\protect{\cite{Berenger1994}} introduced the concept of perfectly matched layer starting from physical considerations on electromagnetic waves. Bérenger modified the Maxwell equations in the absorbing layer so that waves entering the layer are damped out and no reflections arise at the interface. This is why the layer is referred to as perfectly matched.

However, the original approach of Bérenger was based on a splitting technique that led to numerical instabilities at long-time simulations, and could break the hyperbolicity of the system. Then, if the problem is no longer hyperbolic, but only \emph{weakly} hyperbolic, the lower order terms must be treated carefully because some disturbances may arise at later stages of the simulation.

A new construction of PMLs for hyperbolic systems was proposed by Hagstrom \protect{\cite{Hagstrom2003}}. This technique is based on the modal analysis of the governing equations in Laplace--Fourier space to derive the layer model, and is viable only for linear low-order terms. The modal solution inside the layer is constructed so that the eigenfunctions of the problem remain the same across the layer, thereby guaranteeing that no reflections arise at the interface, and the absorbing layer is perfectly matched. Appel\"{o} et al. \protect{\cite{Appelo2006}} later deepened the analysis of this technique and established a solid theoretical foundation.
By using this approach, Gao et al.~\protect{\cite{Gao2011a}} constructed a PML for the BGK equations, and this is the model we study and enhance in the present work.

The reason why we develop and study a PML for the BGK model is that the BGK equations are \emph{linear}, so we can apply the theory of~\protect{\cite{Hagstrom2003,Appelo2006}}. The Navier--Stokes equations (NSEs) are \emph{nonlinear}, so they do not allow the use of such a theory. The ultimate goal would be to couple the absorbing layer developed using the BGK model and the NSEs. Indeed, there exist some formulas that relate the BGK variables to the physical variables modeled by the NSEs. The coupling of the NSE and the BGK equations is left for future work.

\subsection{Contributions}

This work aims to analyze the PML for the BGK equations proposed by Gao et al. \protect{\cite{Gao2011a}} to investigate the role and importance of the parameters appearing in the model. To this aim, we use both analytical and numerical tools.

In particular, the main contributions of this work are as follows:
\begin{itemize}
\item We study, analytically and numerically, the stability of a PML for the BGK model via enforcing the energy decay and a technique based on continued fractions.
\item We apply the ANOVA expansion of multivariate functions to the BGK model with the PML, on two different test cases, and we use it to compute some sensitivity measures to identify the most crucial parameters, their role, and their importance.
\end{itemize}

\subsection{Outline of the paper}

The rest of this work is organized as follows. Section~\ref{sec:BGK_review} presents the Bhatnagar--Gross--Krook (BGK) model of the Boltzmann equation and discusses some implementation aspects. Section~\ref{sec:PML_for_BGK} introduces a PML for the BGK model following \protect{\cite{Gao2011a}}. Section~\ref{sec:stability_analysis} establishes appropriate stability conditions for the BGK model with the PML. Section~\ref{sec:sensitivity_analysis} presents a machinery based on the ANOVA expansion to systematically explore the parameter space and compute the total sensitivity indices. This is done for two different test cases: a wave problem with a localized initial density distribution, and an isentropic vortex problem. Finally, we conclude the paper by summarizing the contributions and providing future research outlooks in Section~\ref{sec:conclusions}.

\section{Review of the BGK model}\label{sec:BGK_review}

In this section, we present the BGK model, its approximation via Hermite polynomials, and show that this results in a symmetric hyperbolic system. The main references for this section are \protect{\cite{Tolke:2000,evans2010partial}}.

The BGK model is an approximation to the Boltzmann equation for rarefied gases. It is named after Bhatnagar, Gross and Krook, who introduced it in their pioneering work~\protect{\cite{BGK:1954}}. The BGK model reads~\protect{\cite[(1)]{Tolke:2000}}\protect{\cite[(6.6, 1)]{Chapman_Cowling:1995}}
\begin{equation}\label{eq:BGK_model}
   \dfrac{\partial f}{\partial t} + \csibold \cdot \nablabold_{\xbold} f = - \dfrac{1}{\tau} \left( f - f_{\mathrm{B}}(\rho,\ubold)\right),
\end{equation}
where $f \equiv f(t,\csibold,\xbold)$ is the particle distribution function, $\csibold = \left[ \xi_{1}, \xi_{2}, \xi_{3} \right]$ is the microscopic velocity, and $\tau$ is a relaxation time. The Maxwell--Boltzmann equilibrium distribution function $f_{\mathrm{B}}$ is defined by~\protect{\cite[(6)]{Tolke:2000}}
\[
   f_{\mathrm{B}}(\rho,\ubold) = \dfrac{\rho}{(2\pi\RT)^{d/2}} \exp \left( - \dfrac{\vert \csibold - \ubold \vert^{2}}{2\RT } \right),
\]
where $\rho$ and $\ubold$ are the macroscopic density and velocity, $R$ is the gas constant, $T$ is the thermodynamic temperature and $d$ is the number of space dimensions. The two terms on the left-hand side of \eqref{eq:BGK_model} represent the mixing and transport of the particles, respectively, while the right-hand side considers the collisions between the particles. 

The relationships between $f(t,\csibold,\xbold)$ and the macroscopic quantities, i.e., density $\rho$, momentum $\rho \ubold $ and pressure tensor $ P_{ij} $, are given by \protect{\cite[(2.2)--(2.6)]{Grad1949}}
\begin{equation}\label{eq:density_velocity}
   \rho = \int_{-\infty}^{+\infty} f \,\mathrm{d}\csibold, \qquad \rho u_{i} = \int_{-\infty}^{+\infty} \zeta_{i}f\,\mathrm{d}\csibold, \qquad P_{ij} = \int_{-\infty}^{+\infty} \left( \zeta_{i} - u_{i} \right) \left( \zeta_{j} - u_{j} \right) f \,\mathrm{d}\csibold.
\end{equation}
The stress tensor $\sigma_{ij}$ is defined by
\begin{equation}\label{eq:stress_tensor}
   \sigma_{ij} = p\,I - P_{ij},
\end{equation}
where $ p = \frac{1}{3}\trace\lbrace P_{ij} \rbrace = \RT\,\rho$ is the scalar pressure.

\subsection{Approximation of the BGK model}

The BGK approximation to the Boltzmann equation is more general than the Navier-–Stokes equations, but it has also more independent variables.
Directly solving the BGK approximation to the Boltzmann equation is too demanding because of the 6+1 dimensions. Hence, as proposed by~\protect{\cite{Tolke:2000}}, we approximate it with an expansion in Hermite polynomials. To recover the macroscopic flow properties, the order of the polynomial space used to approximate the phase space must be sufficiently large. As it was done in~\protect{\cite[\S 2]{Tolke:2000}} and \protect{\cite[\S 2.1]{Karakus:2019}}, to model isothermal and nearly incompressible flows, we use second-order Hermite polynomials. For two-dimensional problems as the ones treated in this paper, a second-order phase space approximation yields a vector of six unknown polynomial coefficients. Therefore, expanding the particle distribution function in a basis of Hermite polynomials $\xi_{k}(\csibold)$~\protect{\cite[(13)]{Tolke:2000}}, we get
\[
   f(t,\csibold,\xbold)= \dfrac{\rho}{(2\pi\RT)^{d/2}} \exp \left( - \dfrac{\csibold \cdot \csibold}{2\RT } \right) \sum_{k=1}^{6} a_{k}(\xbold,t) \, \xi_{k}(\csibold).
\]
After some manipulations, we find the \emph{approximate form of the BGK model}~\protect{\cite[(34)--(39)]{Tolke:2000}}
\begin{align*}
   \frac{\partial a_{1}}{\partial t} + \sqrt{\RT} \left( \frac{\partial a_{2}}{\partial x_{1}} + \frac{\partial a_{3}}{\partial x_{2}} \right) & = 0, \\
   \frac{\partial a_{2}}{\partial t} + \sqrt{\RT} \left( \frac{\partial a_{1}}{\partial x_{1}} + \sqrt{2} \, \frac{\partial a_{5}}{\partial x_{1}} + \frac{\partial a_{4}}{\partial x_{2}} \right) & = 0, \\
   \frac{\partial a_{3}}{\partial t} + \sqrt{\RT} \left( \frac{\partial a_{4}}{\partial x_{1}} + \frac{\partial a_{1}}{\partial x_{2}} + \sqrt{2} \, \frac{\partial a_{6}}{\partial x_{2}}\right) & = 0, \\
      \frac{\partial a_{4}}{\partial t} + \sqrt{\RT} \left( \frac{\partial a_{3}}{\partial x_{1}} + \frac{\partial a_{2}}{\partial x_{2}} \right) & = -\frac{1}{\tau} \left( a_{4} - \frac{a_{2} a_{3}}{a_{1}} \right), \\
      \frac{\partial a_{5}}{\partial t} + \sqrt{2\RT} \, \frac{\partial a_{2}}{\partial x_{1}} & = -\frac{1}{\tau} \left( a_{5} - \frac{a_{2}^{2}}{\sqrt{2} \, a_{1}} \right), \\
      \frac{\partial a_{6}}{\partial t} + \sqrt{2\RT} \, \frac{\partial a_{3}}{\partial x_{2}} & = -\frac{1}{\tau} \left( a_{6} - \frac{a_{3}^{2}}{\sqrt{2} \, a_{1}} \right).
\end{align*}
We can rewrite this system in a compact form as
\begin{equation}\label{eq:BGK_model_approx}
   \dfrac{\partial \a}{\partial t} + A_1 \dfrac{\partial \a}{\partial x_1} + A_2 \dfrac{\partial \a}{\partial x_2} = S(\a),
\end{equation}
where $\a = (a_{1}, a_{2}, a_{3}, a_{4}, a_{5}, a_{6})\tr$ is the vector collecting the expansion coefficients, and $A_{1}$, $A_{2}$ are real symmetric matrices defined by
\[
   A_1 = \sqrt{\RT} 
      \begin{bmatrix}
         \ 0 & 1 & 0 & 0 & 0 & 0 \ \vspace{2mm}\\
         \ 1 & 0 & 0 & 0 & \sqrt{2} & 0 \ \vspace{2mm}\\
         \ 0 & 0 & 0 & 1 & 0 & 0 \ \vspace{2mm}\\
         \ 0 & 0 & 1 & 0 & 0 & 0 \ \vspace{2mm}\\
         \ 0 & \sqrt{2} & 0 & 0 & 0 & 0 \ \vspace{2mm}\\
         \ 0 & 0 & 0 & 0 & 0 & 0 \
      \end{bmatrix}, \qquad 
   A_2 = \sqrt{\RT} 
      \begin{bmatrix}
         \ 0 & 0 & 1 & 0 & 0 & 0 \ \vspace{2mm}\\
         \ 0 & 0 & 0 & 1 & 0 & 0 \ \vspace{2mm}\\
         \ 1 & 0 & 0 & 0 & 0 & \sqrt{2} \ \vspace{2mm}\\
         \ 0 & 1 & 0 & 0 & 0 & 0 \ \vspace{2mm}\\
         \ 0 & 0 & 0 & 0 & 0 & 0 \ \vspace{2mm}\\
         \ 0 & 0 & \sqrt{2} & 0 & 0 & 0 \
      \end{bmatrix}.
\]
The right-hand side term $S(\a)$ in~\eqref{eq:BGK_model_approx} is a nonlinear source vector, defined by
\[
   S(\a) = -\dfrac{1}{\tau} \left( 0, \ 0, \ 0, \ a_{4} - \dfrac{a_{2} a_{3}}{a_{1}}, \ a_{5} - \dfrac{a_{2}^{2}}{\sqrt{2} a_{1}}, \ a_{6} - \dfrac{a_{3}^{2}}{\sqrt{2} a_{1}} \right)\tr.
\]
We can split $S(\a)$ into its linear and nonlinear parts as
\[ 
   S(\a) = S_{\mathrm{L}}(\a) + S_{\mathrm{NL}}(\a),
\]
where
\[
   S_{\mathrm{L}}(\a) = -\dfrac{1}{\tau}\left(0, \ 0,\ 0, \ 1, \ 1, \ 1 \right)\a
\]
collects the linear terms, while
\[
   S_{\mathrm{NL}}(\a) = -\dfrac{1}{\tau} \left( 0, \ 0, \ 0, \ - \dfrac{a_{2} a_{3}}{a_{1}}, \ - \dfrac{a_{2}^{2}}{\sqrt{2} a_{1}}, \ - \dfrac{a_{3}^{2}}{\sqrt{2} a_{1}} \right)\tr
\]
collects the nonlinear terms.
Hence, we can rewrite~\eqref{eq:BGK_model_approx} highlighting the nonlinear terms by moving them to the right-hand side
\begin{equation}\label{eq:BGK_model_approx2}
   \dfrac{\partial \a}{\partial t} + A_1 \dfrac{\partial \a}{\partial x_1} + A_2 \dfrac{\partial \a}{\partial x_2} - S_{\mathrm{L}}(\a) = S_{\mathrm{NL}}(\a).
\end{equation}
In the remaining part of the paper, we may refer to~\eqref{eq:BGK_model_approx} or \eqref{eq:BGK_model_approx2} simply as the BGK model, even though, strictly speaking, the BGK model is~\eqref{eq:BGK_model}, while~\eqref{eq:BGK_model_approx} and \eqref{eq:BGK_model_approx2} are the approximate forms obtained via the Hermite basis expansion.

\subsection{Hyperbolicity}\label{sec:hyperbolicity}

In this section, we look at system \eqref{eq:BGK_model_approx} and investigate its hyperbolic properties following \protect{\cite{Hernquist:1987}} and \protect{\cite[\S 7.3]{evans2010partial}}.
Let $P$ be the differential operator of the system~\eqref{eq:BGK_model_approx}, defined by 
\[
   P(\partial/\partial \xbold) \coloneqq - \left[ A_1 \dfrac{\partial}{\partial x_1} + A_2 \dfrac{\partial}{\partial x_2} \right].
\]
In general, $P$ is also a function of $\xbold \equiv (x_1, x_2)$ and $t$, but here, since we deal with a system of partial differential equations (PDEs) with constant coefficients, it is only a function of $\partial/\partial \xbold\equiv(\partial/\partial x_1,\partial/\partial x_2)$.
We define the \emph{principal part} $P_1$ of the differential system by replacing the partial derivatives $\partial / \partial \xbold$ in $P$ with $\n \equiv (n_1,n_2) \in \R^2$, i.e.,
\begin{equation}\label{eq:same_as_An}
   P_1(\n) \coloneqq - \sqrt{\RT}
   \begin{bmatrix}
      \ 0 & n_1 & n_2 & 0 & 0 & 0 \vspace{2mm}\\
      \ n_1 & 0 & 0 & n_2 & \ \sqrt{2}\,n_1 & 0 \vspace{2mm}\\
      \ n_2 & 0 & 0 & n_1 & 0 & \sqrt{2}\,n_2 \  \vspace{2mm}\\
      \ 0 & n_2 & n_1 & 0 & 0 & 0 \vspace{2mm}\\
      \ 0 & \sqrt{2}\,n_1 & 0 & 0 & 0 & 0 \vspace{2mm}\\
      \ 0 & 0 & \sqrt{2}\,n_2 & 0 & 0 & 0 
   \end{bmatrix}.
\end{equation}
Next, we show some properties of system \eqref{eq:BGK_model_approx}, starting from the following definition.

\begin{definition}[\protect{\cite[\S 7.3.1]{evans2010partial}}]
The system \eqref{eq:BGK_model_approx} is called \emph{hyperbolic} if the matrix $P_1(\n)$ is diagonalizable for each $\n \in \R^{2}, \ t\geq 0 $.
\end{definition}

Since $P_1(\n)$ in~\eqref{eq:same_as_An} is a real symmetric matrix, by the spectral theorem, it is diagonalizable through an orthogonal transformation. Hence, \eqref{eq:BGK_model_approx} is a \emph{hyperbolic} system.
There are also two important special cases.

\begin{definition}[\protect{\cite[\S 7.3.1]{evans2010partial}}]
The system of PDEs \eqref{eq:BGK_model_approx} is \emph{symmetric hyperbolic} if each $A_i$ is a symmetric matrix. Moreover, \eqref{eq:BGK_model_approx} is \emph{strictly hyperbolic} if for any $\n \in \R^2$, $\n \neq 0$ and $t\geq 0 $, all the eigenvalues of $P_1(\n)$ are real and distinct.
\end{definition}

In our case, the matrices $A_1$ and $A_2$ are both symmetric, and hence \eqref{eq:BGK_model_approx} is a \emph{symmetric} hyperbolic system.
To check whether it is also strictly hyperbolic, we look at the eigenvalues of $P_1(\n)$:
\begin{equation}\label{eq:BGK_model_eigenvalues}
   0, \quad 0, \quad  - \sqrt{\RT}\sqrt{n_1^2 + n_2^2}, \quad  \sqrt{\RT}\sqrt{n_1^2 + n_2^2}, \quad  - \sqrt{3\RT}\sqrt{n_1^2 + n_2^2}, \quad  \sqrt{3\RT}\sqrt{n_1^2 + n_2^2}.
\end{equation}
Clearly, \eqref{eq:BGK_model_approx} is not \emph{strictly} hyperbolic since zero appears twice as an eigenvalue of $P_1(\n)$.
The hyperbolicity condition is crucial because it is equivalent to requiring that there are six distinct plane wave solutions of \eqref{eq:BGK_model_approx} for each direction $\n$ \protect{\cite[\S 7.3.1]{evans2010partial}}.

In general, hyperbolic systems may have zero eigenvalues and pairs of eigenvalues with the same magnitude but opposite signs. If we have, say, $p$ pairs of eigenvalues that differ in their sign only, this implies that at each boundary, we can impose at most $p$ boundary conditions. The eigenvalues in \eqref{eq:BGK_model_eigenvalues} tell us that there are six characteristics, among which two are positive, two are negative, and two are zero. According to our previous observation, this means that at any boundary, we do not have to impose more than two boundary conditions. Still, we may impose an additional two corresponding to the non-propagating modes.

\subsection{Relationship with the Navier--Stokes equations}

Using the relationships in \eqref{eq:density_velocity} and \eqref{eq:stress_tensor} and the properties of Hermite polynomials, one can find the following connections between the expansion coefficients $\a$ and the macroscopic quantities~\protect{\cite[(20)--(25)]{Tolke:2000}}
\[
\rho = \int_{-\infty}^{+\infty} f \,\mathrm{d}\csibold = a_{1}, \quad u_{1} = \int_{-\infty}^{+\infty} \xi_{1} f \,\mathrm{d}\csibold = \dfrac{a_{2} \sqrt{\RT}}{a_{1}}, \quad u_{2} = \int_{-\infty}^{+\infty} \xi_{2} f \,\mathrm{d}\csibold = \dfrac{a_{3} \sqrt{\RT}}{a_{1}},
\]
\[
\sigma_{11} = -\int_{-\infty}^{+\infty} \left( \xi_{1} - u_{1} \right)^{2}  f \,\mathrm{d}\csibold + \RT \, \rho = -\RT\left(\sqrt{2} a_{5} - \dfrac{a_{2}^{2}}{a_{1}}\right), \]
\[
\sigma_{22} = -\int_{-\infty}^{+\infty} \left( \xi_{2} - u_{2} \right)^{2}  f \,\mathrm{d}\csibold + \RT \, \rho = -\RT\left(\sqrt{2} a_{6} - \dfrac{a_{3}^{2}}{a_{1}}\right),
\]
\[
\sigma_{12} = -\int_{-\infty}^{+\infty} \left( \xi_{1} - u_{1} \right)\left( \xi_{2} - u_{2} \right)  f \,\mathrm{d}\csibold = -\RT\left(a_{4} - \dfrac{a_{2} a_{3}}{a_{1}}\right),
\]
and the inverse relationships
\begin{equation}\label{eq:a0a1a2}
   a_{1} = \rho, \quad a_{2} = \dfrac{u_{1} \rho}{\sqrt{\RT}}, \quad a_{3} = \dfrac{u_{2}\rho}{\sqrt{\RT}},
\end{equation}
\begin{equation*}
   a_{4} = \dfrac{u_{1} u_{2}\rho - \sigma_{12}}{\RT}, \quad a_{5} = \dfrac{\sqrt{2}}{2}\dfrac{u_{1}^{2}\rho - \sigma_{11}}{\RT},\quad a_{6} = \dfrac{\sqrt{2}}{2}\dfrac{u_{2}^{2}\rho - \sigma_{22}}{\RT}.
\end{equation*}
Following the basic idea from~\protect{\cite{Grad1949}}, it is possible to show that from~\eqref{eq:BGK_model_approx} one can recover the Navier--Stokes equations, under the assumptions that the relaxation time and the Mach number go to zero, namely in the case of weakly compressible flows only~\protect{\cite{Tolke:2000}}.
Let us consider three time scales $\tau$, $\Gamma_0$, $\Gamma_1$ with the relation $\tau \ll\Gamma_0 \ll \Gamma_1$. Here $\tau$ is of the order of magnitude of the collision time, $\Gamma_0$ represents an intermediate time scale, small enough to allow to consider the macroscopic quantities constant in time, and $\Gamma_1$ is the macroscopic time scale on which variations in density and momentum appear. On the scale $\Gamma_0$, under the condition that $\tau$ is very small, the coefficients $(a_{1}, a_{2}, a_{3})$ can be considered constant in time. Moreover, one can get a relation between the stresses and the flow field via a kinematic viscosity $\nu = \RT \, \tau$ and the ideal gas law $ p = \RT \, \rho $. The coefficients $(a_{4}, a_{5}, a_{6})$ are related to the macroscopic variables via~\protect{\cite[(42)--(44)]{Tolke:2000}}
\[
a_{4} = -\tau \left( \dfrac{\partial \rho u_{2}}{\partial x_1} +  \dfrac{\partial \rho u_{1}}{\partial x_2} \right) + \dfrac{u_{1} u_{2}\rho}{\RT},
\]
\begin{equation}\label{eq:a3a4a5}
a_{5} = -\tau \sqrt{2}\,\dfrac{\partial \rho u_{1}}{\partial x_1} + \dfrac{u_{1}^2 \rho}{\sqrt{2}\RT},
\end{equation}
\[
a_{6} = -\tau \sqrt{2}\,\dfrac{\partial \rho u_{2}}{\partial x_2} + \dfrac{u_{2}^2 \rho}{\sqrt{2}\RT},
\]
Substituting~\eqref{eq:a0a1a2} and~\eqref{eq:a3a4a5} into the first three equations of the BGK model, one can recover
\[
\dfrac{\partial \rho}{\partial t} + \dfrac{\partial \rho u_{1}}{\partial x_1} + \dfrac{\partial \rho u_{2}}{\partial x_2} = 0,
\]
\[
\dfrac{\partial \rho u_{1}}{\partial t} + \dfrac{\partial \rho u_{1}^2}{\partial x_1} + \dfrac{\partial \rho u_{1} u_{2}}{\partial x_2} = \dfrac{\partial \sigma_{11}}{\partial x_1} + \dfrac{\partial \sigma_{12}}{\partial x_2} - \dfrac{\partial p}{\partial x_1} ,
\]
\[
\dfrac{\partial \rho u_{2}}{\partial t} + \dfrac{\partial \rho u_{1} u_{2}}{\partial x_1} + \dfrac{\partial \rho u_{2}^2}{\partial x_2} = \dfrac{\partial \sigma_{12}}{\partial x_1} + \dfrac{\partial \sigma_{22}}{\partial x_2} - \dfrac{\partial p}{\partial x_2} ,
\]
with the stress tensor being
\[
\sigma_{ij} = \RT\,\tau \left( \dfrac{\partial \rho u_i}{\partial x_j} + \dfrac{\partial \rho u_j}{\partial x_i} \right).
\]
We recognize the above equations as the two-dimensional isentropic Navier--Stokes equations for a weakly compressible flow.

\section{A PML for the BGK model}\label{sec:PML_for_BGK}

The PML technique was initially introduced by Bérenger \protect{\cite{Berenger1994}}, starting from physical considerations on electromagnetic waves. Bérenger modified Maxwell's equations so that waves getting into the absorbing layer decayed without reflections at the interface. In his original formulation, Bérenger adopted a splitting technique of Maxwell's equations. It was later shown that such a splitting technique breaks the hyperbolicity of the system, leading to numerical instabilities in long-time simulations.

In 2003, Hagstrom \protect{\cite{Hagstrom2003}} proposed a new technique for developing PMLs for hyperbolic systems. This approach, based on the modal analysis in Laplace--Fourier space, enforces the numerical solutions to decay as they propagate into the PML. 
By following in the steps of Hagstrom, Appel\"{o} et al. \protect{\cite{Appelo2006}} developed the fundamental theoretical analysis of PMLs for linear hyperbolic systems, providing general tools to establish stability and well-posedness.

In this section, we follow the approach of \protect{\cite{Hagstrom2003,Appelo2006}}, introduce a PML for the BGK model, verify the matching properties of the PML, study the occurrence of the parameters in the PML, and define the damping functions.

We emphasize that the approach of \protect{\cite{Hagstrom2003,Appelo2006}} is applicable only to \emph{linear} hyperbolic systems and that our BGK model \eqref{eq:BGK_model_approx2} contains a nonlinear term.
Nonetheless, if we neglect the nonlinear term $S_{\mathrm{NL}}(\a)$ appearing in \eqref{eq:BGK_model_approx2}, then we obtain a linear hyperbolic system, and we can apply the technique of \protect{\cite{Hagstrom2003,Appelo2006}} to construct a PML for this problem.

This approach was adopted by~\protect{\cite{Gao2011a}}, who neglected the nonlinear term in \eqref{eq:BGK_model_approx2}, constructed a PML according to~\protect{\cite{Hagstrom2003,Appelo2006}}, and finally appended the nonlinear term at the equations. They proposed the following PML formulation for \eqref{eq:BGK_model_approx} \protect{\cite[(23)]{Gao2011a}}
\begin{equation}\label{eq:BGKPML_full}
   \begin{cases}
      \dfrac{\partial \a }{\partial t} + A_1 \left( \dfrac{\partial \a}{\partial x_1} + \SI \left( \LO \a + \OMEGA\right) \right) + A_2 \left( \dfrac{\partial \a}{\partial x_2} + \SII \left( \LOt \a + \THETA\right) \right) = S(\a),\vspace{3mm}\\
      \dfrac{\partial \OMEGA }{\partial t} + \AI \dfrac{\partial \OMEGA}{\partial x_2} + (\AO + \SI) \OMEGA + \dfrac{\partial \a}{\partial x_1} + \LO (\AO + \SI ) \a - \LI \dfrac{\partial \a}{\partial x_2} = \zerovector,\vspace{3mm}\\
      \dfrac{\partial \THETA }{\partial t} + \AIt \dfrac{\partial \THETA}{\partial x_1} + (\AOt + \SII) \THETA + \dfrac{\partial \a}{\partial x_2} + \LOt (\AOt + \SII ) \a - \LIt \dfrac{\partial \a}{\partial x_1} = \zerovector,
   \end{cases}
\end{equation}
where $\OMEGA$ and $\THETA$ are auxiliary variables and $ \AO$, $\LO$, $\AI$, $\LI$, $\AOt$, $\LOt$, $\AIt$, and $\LIt$ are parameters of the model. The damping functions in the $x$- and $y$-directions are $\SI$ and $\SII$, respectively. 

\subsection{Damping functions}
The positive functions $\SI $ and $ \SII $ appearing in~\eqref{eq:BGKPML_full} are the damping functions, assumed to be smooth and equal to zero only at the PML interface. In general, the damping functions have the form \protect{\cite[(22)]{Gao2011a}}
\begin{equation}\label{eq:damping_function}
   \sigma(x) = C \left( \dfrac{x - x_0}{L}  \right)^{\beta},
\end{equation}
where $x_0$ represents the abscissa at which the PML begins, $L$ is the thickness of the layer, and the exponent $\beta$ is used to control the smoothness of the absorption profile. The constant $C$ represents the overall strength of the absorption and is usually chosen as the inverse of the time-step, $C \approx \left( \Delta t \right)^{-1} $, to avoid restrictions on the time-step caused by the PML.

In the remaining part of this paper, we refer to the BGK model coupled with a PML, equation \eqref{eq:BGKPML_full}, as the \emph{BGK+PML model}.

\subsection{Perfect matching}

The key idea behind the PML is that the eigenfunctions for the eigenvalue problem \emph{inside} the layer have to be the same as \emph{outside} the layer. This is the most straightforward way to design a PML so that reflections at the PML interface are prevented \protect{\cite{Hagstrom2003}}. In the following, we verify this is the case for the BGK+PML model~\eqref{eq:BGKPML_full}.

Consider the homogeneous case of the first equation in \eqref{eq:BGKPML_full} for a PML developing in the $x_1$-direction only, with $\SI$ being a constant for simplicity. The governing equation \emph{outside the layer} is
\[
   \dfrac{\partial \a }{\partial t} + A_1 \dfrac{\partial \a}{\partial x_1} + A_2 \dfrac{\partial \a}{\partial x_2} = \zerovector,
\]
whose Laplace--Fourier transform is
\begin{equation}\label{eq:LFT_BGK}
   \left( s I + \lambda A_1 + \I k_2 A_2 \right) \eigenphi(x_1,\I k_2,s) = \zerovector,
\end{equation}
with modal solution
\[
   \aHAT = e^{\lambda x_1} \eigenphi(x_1,\I k_2,s).
\]
\emph{Inside the layer}, the governing equation is 
\begin{equation}\label{eq:2.3}
\dfrac{\partial \a }{\partial t} + A_1 \left( \dfrac{\partial \a}{\partial x_1} + \SI \left( \LO \a + \OMEGA\right) \right) + A_2 \dfrac{\partial \a}{\partial x_2} = \zerovector.
\end{equation}
This equation has been constructed based on the following \emph{ansatz} for the modal solution inside the layer
\begin{equation}\label{eq:2.4}
   \aHAT_{\mathrm{PML}} = e^{\lambda x_1 + \left[ \frac{\lambda - \LI \I k_2 + \LO \AO}{s + \AI \I k_2 + \AO} - \LO \right] \SI x_1} \eigenphi(x_1,\I k_2,s).
\end{equation}
It can be shown that the Laplace--Fourier transform of \eqref{eq:2.3} is 
\[
   \left( s I + A_1 \left( \left( I - \dfrac{\SI}{\hat{r} + \SI} \right) \left( \dfrac{\partial}{\partial x_1} + \SI \LO\right) + \dfrac{\SI}{\hat{r} + \SI} \left(  \LI \I k_2 - \LO \AO\right)  \right) + \I k_2 A_2 \right) \aHAT_{\mathrm{PML}} = \zerovector.
\]
By inserting the \emph{ansatz}~\eqref{eq:2.4} into the last equation and approaching the PML interface (i.e., letting $\SI \to 0$), one recovers \eqref{eq:LFT_BGK}. This is precisely what we wanted: the eigenfunctions for the governing equations recast in the Laplace--Fourier space remain the same across the PML interface.

In all their simulations, Gao et al. \protect{\cite{Gao2011a}} chose the parameters as
\[
   \LI = 0, \quad \LO = 0, \quad \AI = 0, \quad \AO \neq 0,
\]
\[
   \LIt = 0, \quad \LOt = 0, \quad \AIt = 0, \quad \AOt \neq 0.
\]
The choice of these parameters and the fact that their precise role was not well understood left room for further study that we pursue in this work.

\subsection{Our study case: \texorpdfstring{$\sigma_{2}=0$}{TEXT}}

From now on, we will always consider a PML along the $x_1$-direction only. This is equivalent to set $\SII = 0$, so that the system \eqref{eq:BGKPML_full} becomes
\begin{equation}\label{eq:BGKPML_sigma2=0}
   \begin{cases}
      \dfrac{\partial \a }{\partial t} + A_1 \left( \dfrac{\partial \a}{\partial x_1} + \SI \left( \LO \a + \OMEGA\right) \right) + A_2 \dfrac{\partial \a}{\partial x_2} = S(\a),\vspace{3mm}\\
      \dfrac{\partial \OMEGA }{\partial t} + \AI \dfrac{\partial \OMEGA}{\partial x_2} + (\AO + \SI) \OMEGA + \dfrac{\partial \a}{\partial x_1} + \LO (\AO + \SI ) \a - \LI \dfrac{\partial \a}{\partial x_2} = \zerovector.
   \end{cases}
\end{equation}
Table~\ref{tab:parameters_occurrence} summarizes the parameters of~\eqref{eq:BGKPML_sigma2=0} and their occurrence in the equations. We refer to the first equation in \eqref{eq:BGKPML_sigma2=0} as the $\a$ equation and the second equation as the $\OMEGA$ equation.

\begin{table}[htbp]
   \caption{Occurrence of the parameters in the BGK+PML model \eqref{eq:BGKPML_sigma2=0}.}
   \label{tab:parameters_occurrence}
   \centering
   \begin{tabular}{cp{0.80\textwidth}}
       Parameter     & Occurrence \\ \toprule
       $\LO$         &  Once in the $\a$ equation and once in the $\OMEGA$ equation, in both cases, controlling the behavior of the linear term in $\a$.  We note that, in the $\OMEGA$ equation, $\LO$ appears as a coefficient of $\a$ only if $\AO \neq 0$.\\ \midrule
       $\LI$         &  Once in the $\OMEGA$ equation, as a multiplying coefficient of the derivative of $\a$ with respect to $x_2$. \\ \midrule
       $\AO$         &  Twice in the $\OMEGA$ equation, the first time as a coefficient of the linear term in $\OMEGA$ and the second time as a coefficient of the linear term in $\a$. We note that, in the $\OMEGA$ equation, $\AO$ appears as a coefficient of $\a$ only if $\LO \neq 0$.\\ \midrule
       $\AI$         &  Once in the $\OMEGA$ equation, as a multiplying coefficient of the derivative of $\OMEGA$ with respect to $x_2$. \\ 
   \end{tabular}
\end{table}

\section{Implementation aspects and testing}\label{sec:BGKPML_testing}

This section describes some aspects of implementing and testing the BGK and BGK+PML models.
We employed fourth-order accurate finite differences for the spatial discretization of the first derivatives, and a fourth-order Runge--Kutta method for the time evolution. The Runge--Kutta method is a multi-stage method in which each intermediate stage is in some sense equivalent to a forward Euler method \protect{\cite{Gottlieb2001}}.

To verify that our numerical scheme can simulate the behavior of some simple flows, we tested it using a simple planar Poiseuille flow. This problem was also used a benchmark test case in~\protect{\cite{Tolke:2000} and is excellent because it has an explicit formula for the solution to the Navier--Stokes equations, so we can verify that the BGK model is able to maintain that solution. For more details on the accuracy tests and on the implementation, we refer the reader to~\protect{\cite{Sutti:2015}}.

\subsection{Spatial discretization}

In this work, for simplicity, we adopted a completely regular grid. A scheme using more general unstructured grids, such as the discontinuous Galerkin method, would be less practical. For instance, when comparing two solutions, they would be living on two different grids, and we would have to introduce some mapping to allow comparison. 
Since the primary goal of this work is to study the PML method, we deem a finite difference scheme appropriate for this purpose, at least at an early stage.

\subsection{Time discretization}

For the time discretization, we employ a fourth-order \emph{Runge--Kutta method}:
\begin{align*}
& \kap_{1} = \f\left(t^n, \ \a^n\right),\vspace{3mm}\\
& \kap_{2} = \f\left(t^n + \frac{\Delta t}{2}, \ \a^n + \frac{\Delta t}{2} \kap_{1} \right),\vspace{3mm}\\
& \kap_{3} = \f\left(t^n + \frac{\Delta t}{2}, \ \a^n + \frac{\Delta t}{2} \kap_{2} \right),\vspace{3mm}\\
& \kap_{4} = \f\left(t^n + \Delta t, \ \a^n + \Delta t \kap_{3} \right),\vspace{3mm}\\
& \a^{n+1} = \a^n + \dfrac{\Delta t}{6} \left( \kap_{1} + 2\kap_{2} + 2\kap_{3} + \kap_{4}\right),
\end{align*}
where $\f$ denotes the right-hand side of the BGK equations after spatial discretization, $\Delta t$ is the time-step, and $\a^n$ are the BGK variables computed at time $t^n$.

\subsection{The CFL condition}

In general, for a $d$-dimensional space, the stability condition for a hyperbolic system is expressed as \protect{\cite[(7.17)]{Rezzolla2005}}
\[
   \Delta t \leq C \min\left(\dfrac{h_{x}}{\sqrt{d}\vert u_{2}\vert }\right),
\]
where $C$ is a constant that depends on the method, $i=1,2,\dots,d$ and $ \vert u_{2} \vert = (\sum_{i=1}^{d} v_i^2)^{1/2}$.
In our case, the stable time step becomes
\[
   \Delta t \leq \dfrac{h_{x}}{\sqrt{3} \sqrt{2} \sqrt{2\RT}} = \dfrac{h_{x}}{2 \sqrt{3\RT}},
\]
where $\sqrt{3}$ is the largest eigenvalue, $\sqrt{2}$ is the square root of the number of space dimensions, and $\sqrt{2\RT}$ is the most significant entry of $A_1$ and $A_2$.

The stability condition in two space dimensions can be viewed as an extension of the well-known result in 1D: the numerical domain of dependence of a time-dependent PDE has to contain the physical domain of dependence \protect{\cite{LeVeque2007}}. Courant, Friedrichs, and Lewy wrote a fundamental paper in 1928 that was the first paper on the stability and convergence of finite difference methods for PDEs.

\subsection{Numerical simulations of BGK and BGK+PML models}

In this section, we discuss and report the results of some simulations to assess the accuracy of the BGK+PML model \eqref{eq:BGKPML_sigma2=0} in reproducing the results of the plain BGK model. We first discuss the parameters used in the simulations with the BGK+PML model and then highlight the main differences with respect to the results given by the plain BGK model.

We consider a square domain $ \Omega \coloneqq L_{x} \times L_{y} = [0,1]^{2}$, with an absorbing PML at the right boundary and wall boundary conditions on the rest of the boundary. We employ 20 grid points in each direction, so the mesh size equals $h_{x} \approx 0.0526 $.

As for the damping function \eqref{eq:damping_function}, we choose a PML thickness of $L=0.40$, the exponent $\beta=4$, and the overall absorption strength is calculated as $C = 1/\Delta t$.

For the initial density distribution, we assume a peak located at the center of the domain $\Omega$, having the form
\begin{equation}\label{eq:initial_density}
   a_{1}(x,y,t=0) = 1 + 2\left(p_{\mathrm{in}} - p_{\mathrm{out}}\right) \, \exp\!\left( - \varepsilon \sqrt{\left( x-\tfrac{1}{2}L_{x}\right)^2 + \left( y-\tfrac{1}{2}L_{y}\right)^2 }\right),
\end{equation}
The factor $\varepsilon$ at the exponent is set to $\varepsilon = 10$ to ensure a quick spatial decay of the peak since we want its support outside the PML.

At the initial time, the velocity and the auxiliary variables $\OMEGA$ are set equal to zero everywhere in the domain. In principle, we should evolve the auxiliary variables only inside the PML, but it is easier to solve for them everywhere, and this is what we do here.

For the BGK model without the PML, we keep the same data, the only exception being $L_{x}^{\mathrm{BGK}} \approx 2.5 \times L_{x}$ (the exact value depends on the mesh-size $h_{x}$), since we need to guarantee that any wave propagating in the positive $x$-direction has sufficient space to decay without producing reflections. The additional computational cost of simulating a decaying wave without the PML is evident.
We set the parameters of the PML equal to zero except for $\AO$, which is set to 1, i.e.,
\[
   \LI = 0, \quad \LO = 0, \quad \AI = 0, \quad \AO = 1.
\]
As final time of the simulation, we choose $T = 1.00$, which is sufficient to allow the wave to enter into the layer so that we can observe the PML in action.

Fig.~\ref{fig:BGK_Nx_20} shows the contours of the density distribution $a_{1}$ at different times. Panels on the left correspond to the simulation with the plain BGK model without the PML, while panels on the right refer to the simulation with the BGK+PML model.


\begin{figure}[htbp]
   \centering
   \includegraphics[width=\textwidth]{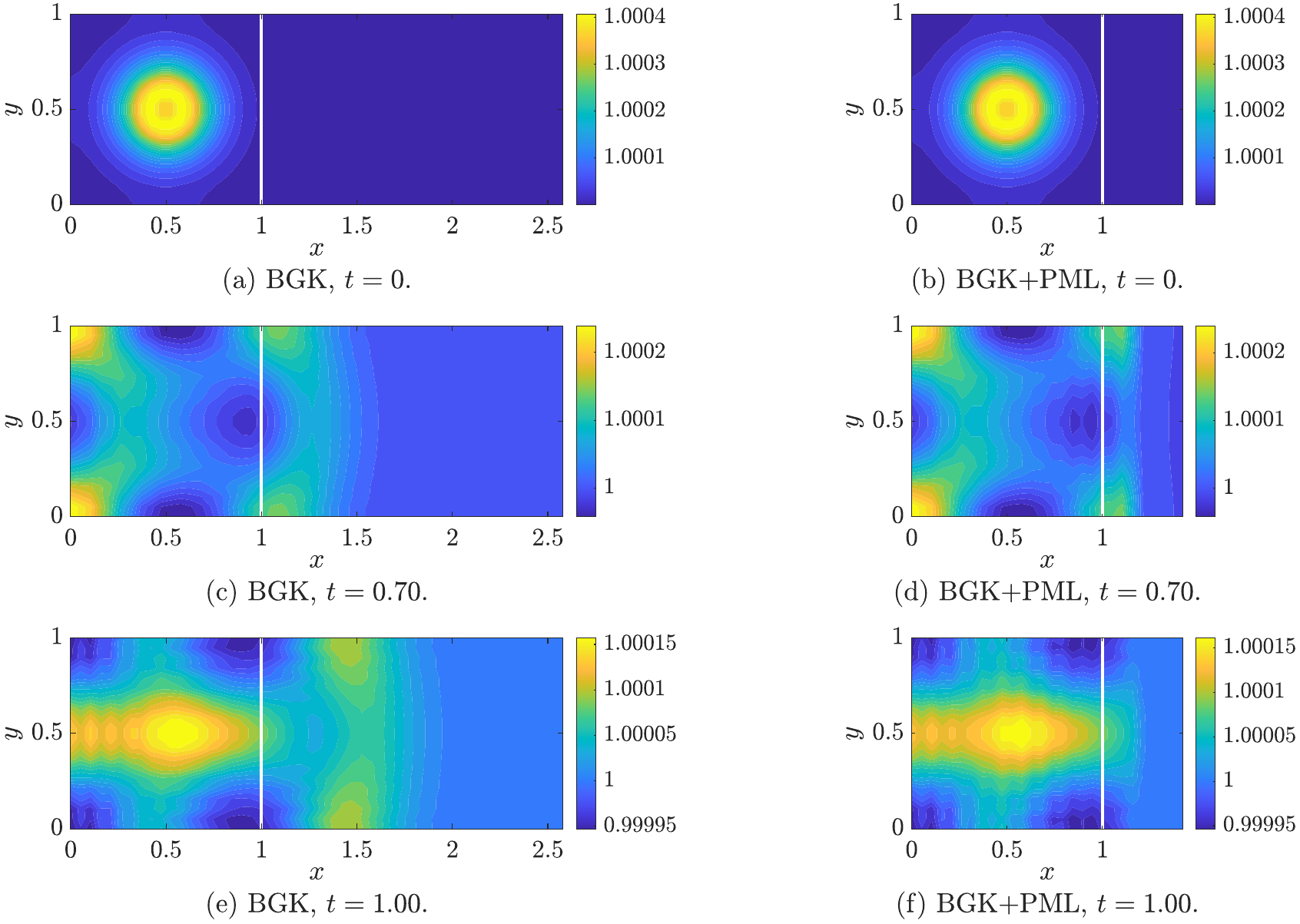}
   \caption{Contour plots of the density distribution at different simulation times.}\label{fig:BGK_Nx_20}
\end{figure}

Panels (a) and (b) of Fig.~\ref{fig:BGK_Nx_20} show the density contours for the initial density distribution. As mentioned above, the peak in the initial density is entirely located in the original domain $\Omega$ outside the PML.
Panels (c) and (d) show the density contours at $t=0.70$. At this time, the waves have already entered the PML, and the latter has already begun to dampen the waves, as can be seen by comparing panels (c) and (d). The waves do not decay immediately as they enter the PML due to the damping function's shape, and the absorption strength becomes more significant as the waves further penetrate the PML.
Panels (e) and (f) show the density contours at the final time of the simulation. Panel (e) highlights the presence of waves propagating to the right, while in panel (f), those waves have been damped out thanks to the PML. Comparing the solution on the domain $\Omega$ in panel (e) with the one in panel (f), it seems that, from a qualitative point of view, the BGK+PML model is well-behaved and is an excellent approximation to the BGK model without the PML. As expected, not only have the waves entering the PML been absorbed, but most importantly, they have not affected the solution on the original domain $\Omega$. We can only notice some minor differences and no significant signs of reflection. The quantitative analyses are carried out in the remaining sections.

As mentioned above, the roles of the parameters in the BGK+PML model are not yet well understood. In the following sections, we will deepen our understanding of the BGK+PML model \eqref{eq:BGKPML_sigma2=0}, first by analyzing stability and establishing reasonable bounds on the parameters, and then by performing a sensitivity analysis in Section~\ref{sec:sensitivity_analysis}.


\section{Stability analysis}\label{sec:stability_analysis}

In this section, we study the stability of the BGK+PML model with two techniques: enforcing the energy decay and via continued fraction expansion.
In both techniques, a vital role is played by the \emph{symbol} of the differential operator of the system.
The analyses provide some reasonable bounds on the parameters of the BGK+PML model. These bounds will then be used in Section~\ref{sec:sensitivity_analysis} to carefully choose the parameters for the simulations needed to perform a sensitivity analysis.

\subsection{The symbol of the BGK+PML model}

The system~\eqref{eq:BGKPML_full} governing the BGK+PML model can be rewritten in matrix form as
\begin{scriptsize}
   \[
      \dfrac{\partial}{\partial t}
         \begin{bmatrix}
            \a \\
            \OMEGA \\
            \THETA
         \end{bmatrix} =
         -\begin{bmatrix}
            A_1 \left( \dfrac{\partial}{\partial x_1} + \SI \LO \right) + A_2 \left( \dfrac{\partial}{\partial x_2} + \SII \LOt \right) & A_1 \SI & A_2 \SII \vspace{3mm}\\
            I\left(\dfrac{\partial}{\partial x_1} + \LO (\AO + \SI ) - \LI \dfrac{\partial}{\partial x_2}\right) & I\left(\AI \dfrac{\partial}{\partial x_2} + \AO + \SI \right) & O \vspace{3mm}\\
            I\left(\dfrac{\partial}{\partial x_2} + \LOt (\AOt + \SII ) - \LIt \dfrac{\partial}{\partial x_1}\right) & O & I\left(\AIt \dfrac{\partial}{\partial x_1} + \AOt + \SII \right)
         \end{bmatrix}
         \begin{bmatrix}
            \a \\
            \OMEGA \\
            \THETA
         \end{bmatrix}.
   \]
\end{scriptsize}
We refer to the matrix on the right-hand side as the \emph{differential operator} of the system, and we denote it by $P \equiv P(\partial/\partial x_1, \partial/\partial x_2)$. Moreover, we denote $[ \a, \ \OMEGA, \ \THETA ]\tr$ by $\ubold \equiv \ubold(x_1,x_2,t)$. Hence, we are dealing with a general system of the type
\begin{equation}\label{eq:PDE_system_short}
   \ubold_{t} = P \ubold,
\end{equation}
with initial condition
\[
   \ubold(x_1,x_2,t=0) = \f(x_1,x_2).
\]
The stability of this problem can be studied via Fourier analysis (sometimes also called Von Neumann analysis). The Fourier transform of \eqref{eq:PDE_system_short} is
\begin{equation}\label{eq:ODE_system}
   \dfrac{\mathrm{d} \uboldHAT}{\mathrm{d} t} = \hat{P} \uboldHAT,
\end{equation}
with initial condition
\[
   \uboldHAT(k_1,k_2,t=0) = \boldsymbol{\widehat{f}}(k_1,k_2),
\]
where $\uboldHAT \equiv \uboldHAT(k_1,k_2,t)$ are the \emph{modes} and $\hat{P}\equiv \hat{P}(\I k_1,\I k_2)$ is called the \emph{symbol} of the differential operator $P$. Equation \eqref{eq:ODE_system} is a recasting of \eqref{eq:PDE_system_short} in the frequency domain, with $k_1,k_2$ being the Fourier variables. Furthermore, \eqref{eq:ODE_system} is a system of ordinary differential equations with constant coefficients, having solution~\protect{\cite[(7.3.38)]{evans2010partial}}, \protect{\cite[(4.5.4)]{Gustafsson1995}}
\begin{equation}\label{eq:mat_exp}
   \uboldHAT = e^{\hat{P}t} \hat{f}(k_1,k_2),
\end{equation}
where $e$ is the matrix exponential. Recasting the original problem \eqref{eq:BGKPML_full} in the frequency domain is useful because the differential problem becomes an \emph{algebraic} problem, and then we can use existing results based on the symbol $\hat{P}$ to establish well-posedness and stability; see, for instance, \protect{\cite[\S 4.5]{Gustafsson1995}}.

Since we are dealing with a wave-dominated problem, it is reasonable to expect that the solution $\uboldHAT$ in \eqref{eq:mat_exp} will decay with time.
In light of this, and because of the properties of the matrix exponential $e^{\hat{P}t}$, the following necessary condition for well-posedness should not be a surprise.

\begin{theorem}[The Petrovskii condition~\protect{\cite[Theorem 4.5.2]{Gustafsson1995}}]\label{thm:petrovskii_condition}
A necessary condition for well-posedness of~\eqref{eq:PDE_system_short} is that, for all $k$, the eigenvalues $\lambda$ of $\hat{P}(\I k)$ satisfy the inequality $ \Re(\lambda) \leq \alpha $, with $\alpha$ being a positive constant.
\end{theorem}

\subsection{Stability analysis via the energy decay}\label{sec:stability_through_energy_decay}

It is possible to work out a stability condition by enforcing energy decay over time. The following stability condition is a standard result.

\begin{theorem}[\protect{\cite[Theorem~4.5.4]{Gustafsson1995}}]
The initial value problem~\eqref{eq:PDE_system_short} is well-posed if there is a constant $\alpha$ such that, for all $k$,
\[
   \hat{P}(\I k) + \hat{P}^{\ast}(\I k) \leq 2\alpha I.
\]
\end{theorem}

We emphasize that this condition is necessary but not sufficient because we \emph{assumed} that the solution is periodic, but in general, it is not. In other words, the above condition is helpful to get a sense of what the parameters do, but it does not give us a complete picture because we have assumed periodicity in space.

The symbol $\hat{P}$ of the differential operator for our problem is
\begin{equation}\label{eq:symbol_diff_op}
   \hat{P} =
   -\begin{bmatrix}
   A_1 (\I k_1 + \SI \LO ) + A_2 (\I k_2 + \SII \LOt ) & A_1 \SI & A_2 \SII \vspace{2mm}\\
   (\I k_1 + \LO (\AO + \SI ) - \I \LI k_2)\,I & (\I \AI k_2 + \AO + \SI )\,I & O \vspace{2mm}\\
   (\I k_2 + \LOt (\AOt + \SII ) - \I \LIt k_1)\,I & O & (\I \AIt k_1 + \AOt + \SII )\,I
   \end{bmatrix}.
\end{equation}
Enforcing the condition $ \hat{P} + \hat{P}^{\ast} < 0 $, and taking into account that $\SI,\SII \geq 0$, we obtain the following conditions on the parameters
\[
   \LO > 0, \quad \LOt > 0, \quad \AO > -\SI, \quad \AOt > -\SII.
\]
We observe that the parameters $\AI$ and $\AIt$ disappear when we take $ \hat{P} + \hat{P}^{\ast} $. Moreover, the parameters $\LI$ and $\LIt$ are involved in the imaginary parts in which $k_1$ and $k_2$ appear; this means that, in principle, they can take any value.

\subsection{Stability analysis via continued fractions}\label{sec:stab_CF}

Appel\"{o} et al.~\protect{\cite{Appelo2006}} proposed an alternative method to analyze stability by studying the sign of the eigenvalues of the symbol $\hat{P}$. In this section, we present the necessary tools for using this method and apply them to~\eqref{eq:symbol_diff_op}. The technique of~\protect{\cite{Appelo2006}} is based on a theorem by Frank~\protect{\cite{Frank1946}}, which is reported here for the reader's convenience. This technique is related to the Routh--Hurwitz stability criterion, often used in control theory. 

\begin{theorem}[Frank, 1946 \protect{\cite[Corollary (38,1b)]{Marden1966}}]\label{thm:frank}
Consider any polynomial $q(z)$ of degree $n$. Let $D$ be a real number and define the polynomials $Q_0$ and $Q_1$ with real coefficients by
\[ q(\I D) \equiv \I^n [ Q_0(D) + \I Q_1(D)] . \]
Then, there is a continued fraction
\[ \dfrac{Q_1(D)}{Q_0(D)} = \dfrac{1}{c_1 D + d_1 - \dfrac{1}{c_2 D + d_2 - \dfrac{1}{c_3 D + d_3 - \cdots - \dfrac{1}{c_{n_{r}} D + d_{n_{r}}}}}}\]
with $c_j \neq 0 $ and $n_{r} \leq n$. The number of roots of $q(z)$ with positive (negative) real part equals the number of positive (negative) $c_j$. Moreover, there are $n-n_{r}$ roots on the imaginary axis.
\end{theorem}

There are several key points to notice about this theorem. First, we must be able to write any polynomial $q(z)$ to read off the polynomials $Q_0$ and $Q_1$ of a real variable. This rewriting can always be achieved. Secondly, we must be able to write the rational function $Q_1/Q_0$ in a continued fraction form. This can always be accomplished via an algorithm that recursively uses the Euclidean division between polynomials, no matter how complicated the rational function of departure, and returns the continued fraction expansion.
Nonetheless, the calculations are far more tedious. Notice also that the total number of coefficients $c_j$ appearing in the above expression may be less than the polynomial degree ($n_{r} \leq n$). Finally, the theorem also tells us that $n - n_r$ roots lie on the imaginary axis.

Without this theorem, it would be too complicated to find an analytical expression for the eigenvalues and study their signs directly. Indeed, suppose we apply Theorem~\ref{thm:frank} to the characteristic polynomial of the symbol $\hat{P}$. In that case, we can determine the sign of its eigenvalues without knowing their exact, explicit form. This is all the information we need to apply the Petrovskii condition, Theorem~\ref{thm:petrovskii_condition}, which requires that all the coefficients $c_j$ in Theorem~\ref{thm:frank} must be defined and negative.

\paragraph{Characteristic polynomial of $\hat{P}$.}
The characteristic polynomial $p(z)$ of the symbol $\hat{P}$~\eqref{eq:symbol_diff_op} factorizes as
\begin{equation}\label{eq:char_pol}
   p(z) = z^{2} \left( z + \AOt + \I k_1 \AIt \right)^{6} \left( z + \AO + \SI + \I k_2 \AI \right)^{2} \mu_{4}(z) \, \nu_{4}(z),
\end{equation}
where $\mu_{4}(z)$ and $\nu_{4}(z)$ are two fourth degree polynomials, defined by
\begin{align*}
   \mu_{4}(z) = & \left( z^2 + k_{1}^{2} + k_{2}^{2} \right) \left( z + \AO + \I k_2 \AI \right)^{2} + 2\left( z + \AO + \I k_2 \AI \right) \times  \\
   & \times \left( k_{2}^{2} + z \left( z - \I k_1 \LO \right) + k_1 k_2 \left( \AI \LO + \LI \right) \right) \SI + \\
   &  + \left( -z^2 \left(-1 + \LO^{2} \right) - 2 \I z k_2 \LO \left( \AI \LO + \LI \right) + k_{2}^{2} \left( 1 + \left( \AI \LO + \LI \right)^{2}\right)\right) \SI^{2},
\end{align*}
\begin{align*}
   \nu_{4}(z) = & \ \left( z^2 + 3 k_{1}^{2} + 3 k_{2}^{2} \right) \left( z + \AO + \I k_2 \AI \right)^{2} + 2\left( z + \AO + \I k_2 \AI \right) \times  \\
   & \times \left(  z \left( z - 3 \I k_1 \LO \right) + 3 k_2 \left( k_2 + k_1 \left( \AI \LO + \LI \right) \right) \right) \SI + \\
   & + \left( z^2 \left(1 -3\LO^{2} \right) - 6 \I z k_2 \LO \left( \AI \LO + \LI \right) + 3 k_{2}^{2} \left( 1 + \left( \AI \LO + \LI \right)^{2}\right)\right) \SI^{2} .
\end{align*}
We observe that the parameters $\LOt$ and $\LIt$ do not appear in the characteristic polynomial~\eqref{eq:char_pol}, which indicates that they can take any value.
From~\eqref{eq:char_pol}, it is clear that $\hat{P}$ has two zero eigenvalues, six times the eigenvalue $z = -(\AOt + \I k_1 \AIt)$, 
and twice the eigenvalue $z = - (\AO + \SI + \I k_2 \AI)$. Since we require the real parts of the eigenvalues to be strictly positive, we obtain the conditions
\[
   \AOt > 0, \qquad \AO > -\SI.
\]
These two conditions coincide with those found in Section~\ref{sec:stability_through_energy_decay} via the energy decay approach.

The polynomials $\mu_{4}(z)$ and $\nu_{4}(z)$ are fourth-degree polynomials, so, in principle, they admit an algebraic solution in closed form. However, the closed-form expression of the solution is too complicated to allow analysis.

In what follows, we apply Theorem~\ref{thm:frank} to $\mu_{4}(z)$ and $\nu_{4}(z)$ separately to work out their respective continued fraction expansions and study the sign of the eigenvalues.

\subsubsection{Application of Theorem \ref{thm:frank} to \texorpdfstring{$ \mu_{4}(z) $}{TEXT}} \label{subsec:mu4_CF}

The first coefficient in the continued fraction expansion of $\mu_{4}(z)$ is
\[
   c_1 = -\dfrac{1}{2(\AO + \SI)}.
\]
It is defined if $\AO \neq -\SI$ and negative if $ \AO > -\SI $, which again coincides with the condition found via the energy decay technique.

The second coefficient in the expansion is
\begin{align*}
   c_2 =  -2(\AO + \SI)^3 / & \left[ \, \AO^{4}+\AO \left( k_{1}^{2} + 4 \AO^{2} - k_1 k_2 ( 2\AI \LO + \LI ) \right) \SI + \right. \\
   & - \left. \left( \AO^{2} (-6 + \LO^{2}) + k_1^{2} (-1 + \LO^{2}) + k_1 k_2 (2\AI \LO + \LI)\right) \SI^{2} + \right. \\
   & - 2\AO (-2+\LO^{2})\,\SI^{3} - \left. (-1+\LO^{2})\,\SI^4 \, \right].
\end{align*}
Due to the complicated expression at the denominator, we analyze this coefficient in two limit cases.

\paragraph{Case $\SI \to 0$.} We first consider the limit case in which we approach the PML interface. In this case, we may drop the higher-order terms in $\SI$, and $c_2$ becomes
\begin{equation*}
   c_2 =  -\dfrac{2(\AO + \SI)^3}{\AO^{4}+\AO \left( k_{1}^{2} + 4 \AO^{2} - k_1 k_2 ( 2\AI \LO + \LI ) \right) \SI}.
\end{equation*}
Then the questions remain: \emph{is $c_2$ defined?} and, if it is defined, \emph{is it negative?}
To answer both these questions, we study when the denominator of $c_2$
\begin{equation}\label{eq:den_c_2}
   f(k_1,k_2) = \AO^{4}+\AO \left( k_{1}^{2} + 4 \AO^{2} - k_1 k_2 ( 2\AI \LO + \LI ) \right) \SI
\end{equation}
is strictly positive. We seek an analytical expression for the boundaries of the instability regions. Starting from $f(k_1,k_2)>0$, we work out $k_2$ as a function of $k_1$ (with $\AO \neq 0$)
\[
   \AO^{3}+ \left( k_{1}^{2} + 4 \AO^{2} - k_1 k_2 ( 2\AI \LO + \LI ) \right) \SI > 0,
\]
\[
   \AO^{3}+ k_{1}^{2} \SI + 4 \AO^{2} \SI > k_1 k_2 ( 2\AI \LO + \LI ) \SI.
\]
If $k_1 = 0$, one has
\[
   \AO^3 + 4 \AO^2 \SI > 0,
\]
which is always guaranteed because $ \AO > -\SI $ (the condition on the first coefficient).
In contrast, for $k_1 \neq 0$, one obtains
\[
   \begin{cases}
      k_2 < \dfrac{\AO^{3}+ \left(4 \AO^{2} + k_{1}^{2} \right) \SI }{ \left( 2\AI \LO + \LI \right) k_1\SI } \quad \mathrm{if} \ k_1>0,\vspace{3mm}\\
      k_2 > \dfrac{\AO^{3}+ \left(4 \AO^{2} + k_{1}^{2} \right) \SI }{ \left( 2\AI \LO + \LI \right) k_1\SI } \quad \mathrm{if} \ k_1<0.
   \end{cases}
\]
Since $k_2$ should be allowed to take any value, the only possibility is that the expression on the right-hand side of the last inequalities is unbounded
\[
   \dfrac{\AO^{3}+ \left(4 \AO^{2} + k_{1}^{2} \right) \SI }{ \left( 2\AI \LO + \LI \right) k_1\SI } \to \infty, \quad \forall k_1 \neq 0,
\]
and since $k_1\SI \neq 0$, we need
\[
   \left( 2\AI \LO + \LI \right) \to 0,
\]
which means that either $\LO = \LI = 0$, or $\AI = \LI = 0$ or $\LO = - \LI/2\AI $. If we assume that all the parameters are positive, then the third option has to be discarded. Moreover, all the preliminary numerical simulations have shown that in practice, $\LO$ has to stay zero to guarantee stability, so we are left only with the first condition
\[
   \LO = \LI = 0.
\]
In general, we observe that the presence of the instability region is associated with the mixed term in $k_1 k_2$ in \eqref{eq:den_c_2}. If somehow the leading coefficient of this mixed term is zero, then $f(k_1,k_2)$ is always positive, and hence $c_2$ is negative.

One can find the same conditions by following the same procedure for $\nu_{4}(z)$. However, the coefficients $c_3$ and $c_4$ of the continued fraction expansions of both $\mu_{4}(z)$ and $\nu_{4}(z)$ have more complicated expressions that do not allow analysis.

Finally, we recall that we only considered the limit case $\SI \to 0$ in which we approach the PML interface. Therefore, the conditions we found above may not give us the complete picture. 
However, the stability analysis is confirmed by the fact that in all the numerical simulations performed with nonzero $\LO$, $\LI$, $\AO$, and $\AI$, the BGK+PML model is unstable. 
In contrast, the simulations performed with $\LO$, $\LI$, $\AO$, and $\AI$ set to zero demonstrated that the BGK+PML model is stable, but one must address accuracy issues. In the next section, we will further develop these aspects.

\section{Sensitivity analysis}\label{sec:sensitivity_analysis}

We are now able to assess how the parameters affect the outcome of the simulations of the BGK+PML model.

We can find many methods in the literature to sample a parameter space and perform a sensitivity analysis. A good review is given in \protect{\cite{Andres1997}}, where the author points out how models with many parameters often behave as if they depend on only a few.
Besides providing accurate results, a suitable sampling method should minimize the number of simulation runs.
The most straightforward way to explore a parameter space would be to generate a simple random sample of the parameters in some given intervals of variation. Unfortunately, random sampling has a convergence rate of $1/\sqrt{d}$, where $d$ is the dimension of parameter space, because of the law of large numbers. Random sampling may prove inefficient and computationally expensive if one has many parameters. Hence, we must find a better way to explore the parameter space systematically. 

The tool we use to gain insights into the role of the parameters is the Analysis of Variance (ANOVA) expansion of multivariate functions and the related concept of total sensitivity indices (TSIs). ANOVA expansions are very useful when one wants to study functionals of solutions to nonlinear partial differential equations. 

In this section, we first describe the ANOVA expansion of multivariate functions and the related concept of TSIs. Then, we define a functional of the solution to the BGK+PML model. We establish some reasonable bounds on the parameters of our BGK+PML model that are involved in the functional chosen. We then apply the ANOVA expansion to our chosen functional and compute the related TSIs. Finally, we explore other possible choices of functional.

\subsection{ANOVA expansion of multivariate functions}

The main references for this section are~\protect{\cite{Saltelli2000,Cao2009,Gao2011}}.
Let $ \Pset =\lbrace 1, 2, \dots, p \rbrace$ be the set of coordinate indices of a $p$-dimensional function, and let $\alphabold = (\alpha_{1}, \alpha_{2}, \dots, \alpha_{p} ) \in \R^p $ be a $p$-dimensional vector. Let $ \Tset \subseteq \Pset $ be a subset of $ \Pset $, and let $t$ denote the cardinality of $ \Tset $ (so $t \leq p$). We denote by $ \alphabold_{\Tset} \in \R^t $ the $t$-dimensional vector that contains the components of $\alphabold \in \R^p $ indexed by $\Tset$. Furthermore, we denote by $A^p$ the $p$-dimensional unit hypercube $[0,1]^p$. 
Any $p$-dimensional function $ g \in L^2(A^p) $ can be written as the \emph{ANOVA expansion} 
\begin{equation}\label{eq:ANOVA_expansion}
   g(\alphabold ) = g_0 + \sum_{\Tset \subseteq \Pset } g_{\Tset}(\alphabold_{\Tset}),
\end{equation}
where the terms in the expansion are calculated recursively through the formula
\begin{equation}\label{eq:ANOVA_recursion}
   g_{\Tset}(\alphabold_{\Tset}) = \int_{A^{p-t}} g_{\Tset}(\alphabold_{\Tset}) \,\mathrm{d}\alphabold_{\Pset\setminus\Tset} - \sum_{\Wset\subset\Tset} g_{\Wset}(\alphabold_{\Wset}) - g_0,
\end{equation}
starting with the zeroth order term
\[
   g_0 = \int_{A^p} g(\alphabold) \,\mathrm{d}\alphabold,
\]
and where, by convention,
\[
   \int_{A^0} g(\alphabold) \,\mathrm{d}\alphabold_{\varnothing} = g(\alphabold).
\]
Each term $g_{\Tset}(\alphabold_{\Tset}) $ in the ANOVA expansion is, in general, a nonlinear function of its $t$ arguments, and it is the \emph{unique} term in the expansion that depends on the $t$ variables indexed by $\Tset$. In other words, the term $g_{\Tset}(\alphabold_{\Tset}) $ describes the effect within $g(\alphabold)$ when those $t$ arguments are simultaneously taken into account.

We emphasize that $\mathrm{d}\alphabold_{\Pset\setminus\Tset}$ in \eqref{eq:ANOVA_recursion} indicates integration over all those coordinate indices \emph{not} included in $\Tset$, and that the sum is carried out over \emph{strict} subsets $ \Wset $ of $ \Tset $. The operation \eqref{eq:ANOVA_recursion} can be regarded as a \emph{projection} since the resulting function depends only on the coordinate indices in $\Tset$. The total number of terms in the ANOVA expansion is $2^p$. 

We stress that the ANOVA expansion is \emph{exact} and contains a \emph{finite} number of terms. If we \emph{truncate} it, we obtain an approximation to $g(\alphabold)$, having fewer terms than the whole expansion. A natural question is how to truncate the ANOVA expansion to get a good approximation to $g(\alphabold)$. This question leads us to the concept of \emph{effective dimension of multivariate functions}; see Sections \ref{subsec:truncated_ANOVA} and \ref{subsec:effective_dimension}.

To clarify equations~\eqref{eq:ANOVA_expansion} and \eqref{eq:ANOVA_recursion}, we explicitly write the expressions to calculate the first few terms in the ANOVA expansion. We define the \emph{order} of a term $g_{\Tset}(\alphabold_{\Tset})$ appearing in~\eqref{eq:ANOVA_expansion} as the cardinality $t$ of the corresponding set $\Tset$. The case $t=1$ generates the \emph{first-order terms}, or \emph{univariate functions}, given by
\[
   g_{i_{1}}(\alpha_{i_{1}}) = \int_{A^{p-1}} g(\alpha) \,\mathrm{d}\alphabold ' - g_0, \qquad i_{1} = 1,2,\dots,p,
\]
where $\mathrm{d}\alphabold '$ indicates integration over all coordinates except $\alpha_{i_{1}}$.

The case $t=2$ generates the \emph{second-order terms}, or \emph{bivariate functions}, given by
\[
   g_{i_{1}i_{2}}(\alpha_{i_{1}},\alpha_{i_{2}}) = \int_{A^{p-2}} g(\alphabold) \,\mathrm{d}\alphabold '' - g_{i_{1}}(\alpha_{i_{1}}) - g_{i_{2}}(\alpha_{i_{2}}) - g_0, \qquad i_{1} < i_{2}, \quad i_{1},i_{2} = 1,2,\dots,p,
\]
where $\mathrm{d}\alphabold ''$ indicates integration over all coordinates except $\alpha_{i_{1}}$ and $\alpha_{i_{2}}$.

The case $t=3$ generates the \emph{third-order terms}, or \emph{trivariate functions}, given by
\begin{equation*}
\begin{split}
   g_{i_{1}i_{2}i_{3}}(\alpha_{i_{1}},\alpha_{i_{2}},\alpha_{i_{3}}) = & \int_{A^{p-3}} g(\alphabold) \,\mathrm{d}\alphabold ''' - g_{i_{1}i_{2}}(\alpha_{i_{1}},\alpha_{i_{2}}) - g_{i_{1}i_{3}}(\alpha_{i_{1}},\alpha_{i_{3}}) - g_{i_{2}i_{3}}(\alpha_{i_{2}},\alpha_{i_{3}}) \\
   & - g_{i_{1}}(\alpha_{i_{1}}) - g_{i_{2}}(\alpha_{i_{2}}) - g_{i_{3}}(\alpha_{i_{3}}) - g_{0}, \\
   & i_{1} < i_{2} < i_{3}, \quad i_{1},i_{2},i_{3} = 1,2,\dots,p,
\end{split}
\end{equation*}
where $\mathrm{d}\alphabold '''$ indicates integration over all coordinates except $\alpha_{i_{1}}$, $\alpha_{i_{2}}$ and $\alpha_{i_{3}}$, and so on. Note that, as we go to higher order, i.e., as $t$ increases, the dimensionality of the integrals that we need to compute decreases. The binomial coefficient gives the total number of $t$th order terms
\[
   \binom{p}{t} = \frac{p!}{t!\,(p-t)!}.
\]
The ANOVA expansion of $g(\alphabold)$ is finally written as
\[
   g(\alphabold) = g_0 + \sum_{i_{1}}^{p} g_{i_{1}} + \sum_{i_{1},i_{2}}^{\binom{p}{2}} g_{i_{1}i_{2}} + \sum_{i_{1},i_{2},i_{3}}^{\binom{p}{3}} g_{i_{1}i_{2}i_{3}} + \dots + \sum_{i_{1},\dots, i_{p-1}}^{p} g_{i_{1}\cdots i_{p-1}} + g_{i_{1}i_{2} \cdots i_{p-1} i_{p}}.
\]
From a computational point of view, the bottleneck in this procedure is represented by the evaluation of the multidimensional integrals needed to construct the expansion. This aspect should not be underestimated because it might compromise the algorithm's efficiency; some aspects in this regard are discussed in Section~\ref{sec:multivariate_numerical_integration}.

\subsubsection{Properties of the ANOVA expansion}

This section provides a list of the most important properties of the ANOVA expansion; a complete list can be found in \protect{\cite{Saltelli2000}}.
\begin{itemize}
\item The ANOVA expansion of a general $p$-dimensional function $g \in L^{2}(A^p)$ is \emph{exact} and \emph{finite}, and contains a total number of $2^p$ terms;
\item the \emph{zero-th order term} $g_0$ in \eqref{eq:ANOVA_expansion} is an integral average of $g$ over the entire parameter space $A^{p}$, and it is a \emph{constant};
\item the term $g_{\Tset}(\alphabold_{\Tset})$ is a function \emph{only} of the coordinates indexed by $\Tset$;
\item \emph{the terms in the ANOVA expansion are mutually orthogonal}, i.e.,
\[
   \int_{A^{p}} g_{\Sset}(\alphabold_{\Sset}) \, g_{\Tset}(\alphabold_{\Tset}) \,\mathrm{d}\alphabold = 0,
\]
whenever $\Sset$ and $\Tset$ contain at least one different index.
This holds also when $\Sset$ and $\Tset$ have the same cardinality.
Note that when $g_{\Sset} = g_0$ we get the particular case
\[
   \int_{A^{p}} g_0 \, g_{\Tset}(\alphabold_{\Tset}) \,\mathrm{d}\alphabold = 0,
\]
which, since $g_0$ is constant, implies
\[
   \int_{A^{p}} g_{\Tset}(\alphabold_{\Tset}) \,\mathrm{d}\alphabold = 0,
\]
meaning that \emph{the terms $g_{\Tset}(\alphabold_{\Tset})$ in the ANOVA expansion have zero average}.
\item each term $g_{\Tset}(\alphabold_{\Tset})$ in the expansion is a \emph{projection} of $g(\alphabold)$ onto a subspace of $L^{2}(A^p)$, with respect to the $L^{2}(A^p)$ inner product. 
\end{itemize}

\subsubsection{The truncated ANOVA expansion} \label{subsec:truncated_ANOVA}

A truncated ANOVA expansion of order $r$ is given by~\protect{\cite[(7)]{Cao2009}}
\[
   g(\alphabold;r) = g_0 + \sum_{\Tset \subseteq \Pset, \ t \leq r } g_{\Tset}(\alphabold_{\Tset}).
\]
The truncated ANOVA expansion $g(\alphabold;r)$ with $r \ll p $ often gives an excellent approximation to $g(\alphabold)$.
Some aspects of well approximating a multivariate function $g(\alphabold)$ by a truncated ANOVA expansion $g(\alphabold;r) $ are presented in \protect{\cite[\S 2.1]{Cao2009}}. For instance, if $r \ll p$, then our $g(\alphabold)$, a function of $p$ arguments, can be well described by a sum of terms, each of which depends at most on $r$ variables. This means that \emph{the contributions provided by coordinate sets having more than $r$ variables can be disregarded}. This leads us to the concept of effective dimension of a function.

\subsubsection{The effective dimension of a function} \label{subsec:effective_dimension}

As mentioned above, the ANOVA expansion is related to the concept of effective dimension of a multivariate function \protect{\cite[\S 2.2]{Cao2009}}.
We define the terms $V_{\Tset}(g)$ and $V(g) $ as
\begin{equation}\label{eq:VTV}
   V_{\Tset}(g) = \int_{A^p} \left( g_{\Tset}(\alphabold_{\Tset}) \right)^2 \,\mathrm{d}\alphabold, \qquad V(g) = \sum_{t>0} V_{\Tset}(g).
\end{equation}
Note that $V_{\Tset}(g)$ is the integral average of the square of the terms appearing in the ANOVA expansion and can be regarded as a variability of $g$ over a given set $\Tset$.

\begin{definition}
The effective dimension of a multivariate function $g$ in the superposition sense or, in short, the \emph{superposition dimension}, is the smallest integer $r$ such that
\[
   \sum_{0<t\leq r} V_{\Tset}(g) \geq q\,V(g),
\]
where $q > 0 $ is called \emph{proportion} and it is typically chosen to be slightly less than $1$; $q = 0.99$ is a common choice.
\end{definition}
\begin{definition}
Given a function $g$ and its approximation $h$, the \emph{normalized approximation error} is defined by
\[
   E(g,h) = \dfrac{1}{V(g)} \int_{A^p} \left( g(\alphabold) - h(\alphabold) \right)^2 \,\mathrm{d}\alphabold.
\]
\end{definition}
We have the following remarkable theorem about the approximation property of the truncated ANOVA expansions.

\begin{theorem}[\protect{\cite[Theorem~2]{Wang:2003}}]
Assume that $g(\alphabold)$ has superposition dimension $r$ in proportion $q$ and let $g(\alphabold;r) = \sum_{0<t<r}g_{\Tset}(\alphabold_{\Tset}) $ denote its truncated ANOVA expansion of order $r$. Then
\[
   E\left( g(\alphabold), g(\alphabold;r)\right) \leq (1 - q).
\]
\end{theorem}

This theorem formalizes what we claimed in the previous section, i.e., if the superposition dimension is small ($r \ll p$), then $g(\alphabold)$ can be well approximated by a truncated ANOVA expansion with only a few terms. More precisely, the error between the ANOVA expansion and its order $r$ truncation is at most $(1-q)$.

Many practical applications have demonstrated that truncated ANOVA expansions of order two can already yield excellent approximations to the original function $g$ \protect{\cite{Cao2009,Gao2011}}. This fact may be explained by the observation that, typically, high-dimensional functions are not really high-dimensional. In practice, multidimensional functions that truly depend on all the parameters are found quite seldom. Usually, the bivariate terms in the ANOVA expansion still matter, but if we consider terms of order higher than two, one can observe that they make a slight difference. The ANOVA expansion is very useful because it quantifies how much structure hides behind a multivariate function.

\subsubsection{Total sensitivity indices}

The total sensitivity index (TSI) of a parameter $\alpha_i$ measures the combined sensitivity of all terms that depend on $\alpha_i$, $i=1,\dots,p$. We define the \emph{sensitivity measure} as
\[
   S_{\Tset} = \dfrac{V_{\Tset}}{V},
\]
where $V_{\Tset}$ and $V$ are defined according to~\eqref{eq:VTV}.
The following result holds
\[
   \sum_{i \in \Tset} S_{\Tset} + \sum_{i \notin \Tset} S_{\Tset} = 1,
\]
where the first term is a sum of the sensitivity measures $S_{\Tset}$ containing the coordinate index $i$, and the second term is a sum of those $S_{\Tset}$ that do not contain it.
We call the first term in the above expression the TSI($i$) of variable $\alpha_i$,
\[
   \text{TSI}(i) \coloneqq \sum_{i \in \Tset} S_{\Tset}.
\]
The TSIs give us a feeling of which parameters are most important. They represent the final goal of our computation of the ANOVA expansion.

\subsection{Multivariate numerical integration}\label{sec:multivariate_numerical_integration}

A crucial aspect in computing the ANOVA expansion is the accurate calculation of the multidimensional integrals appearing in \eqref{eq:ANOVA_recursion}.
The ANOVA expansion theory is based on the assumption that these multidimensional integrals can be evaluated \emph{exactly}. In our case, this is not possible since we do not have an analytic expression for $g$. In practice, we have to resort to multivariate numerical integration.
 The efficiency and the accuracy of the integration methods adopted directly affect the efficiency of the calculations and the accuracy of the ANOVA expansion.

There are several approaches to evaluate multidimensional integrals, for instance,
\begin{itemize}
\item the \emph{Stroud cubature}, which is the simplest approach and gives the minimum amount of nodes to obtain a certain accuracy in high dimensions, but cannot provide very high accuracy;
\item \emph{product rules}, that allow the extension of many known univariate integration formulas to higher dimensions. This approach allows the calculation of the integrals accurately, but it quickly becomes computationally expensive because the number of samples grows like $n^p$ for a quadrature using $n$ points in $p$ dimensions. For instance, the Cartesian product with a Gaussian quadrature for seven dimensions with five integration abscissas in each dimension needs $5^7=78\,125$ evaluations.
\item the Smolyak construction, which is a sparse grid integration method.
\end{itemize}

As our final goal is to compute the TSIs in order to get a sense of which parameters are most significant, the evaluation of the multivariate integrals does not need to be very accurate. We are not interested in the exact or highly accurate values of the TSIs but in understanding how the parameters relate to each other and which one is more important.
In this work, we adopted the Gauss--Legendre quadrature rules, denoted by $(G_{n})^{d}$, where $n$ is the number of nodes, and $d$ the dimension of the function.

\subsection{Definition of the error functional}

We have to define an outcome of the solution on which to focus our analysis.
Preferably, we would choose a functional of the numerical solution to our BGK+PML model. In particular, we choose the maximum over time of the $L^2$-norm of the error in the density $a_{1}$ between the BGK+PML and the plain BGK, calculated along a line close to the PML, normalized with respect to the $L^2$-norm of the initial condition of the density $a_{1}$ on that same line for the plain BGK. Formally,
\begin{equation}\label{eq:error_functional}
   g_{1} = \dfrac{\displaystyle\max_{t\in \left[ 0,T \right]} \left\lbrace \left[ \int_{0}^{L_{y}} \left( a_{1}^{\mathrm{PML}}(x^{\ast},y,t) - a_{1}(x^{\ast},y,t)\right)^2 \dy \right]^{1/2}\right\rbrace}{\left[ \int_{0}^{L_{y}} \left( a_{1}(x^{\ast},y,t=0)\right)^2 \dy \right]^{1/2} },
\end{equation}
where $ T $ is the final time of the simulation, $ x^{\ast} $ is the abscissa at which the reference line is located, while $a_{1}^{\mathrm{PML}}$ and  $a_{1}$ are the densities with and without the PML, respectively. In the most general case we will consider, this error functional will depend on four parameters, namely $\AO$, $\AI$, the PML exponent $\beta$, and the PML thickness $L$. 
In Section~\ref{sec:other_functionals}, we will also consider other functionals, but until then, we will stick to \eqref{eq:error_functional}.

\subsection{Bounds on \texorpdfstring{$\LO$}{TEXT}, \texorpdfstring{$\LI$}{TEXT}, and the PML thickness} \label{sec:LO_LI_to_zero}

In Section~\ref{subsec:mu4_CF}, we have seen from an analytical point of view that to guarantee stability, we must set the parameters $\LO$ and $\LI$ to zero.
Here, we show that the numerical simulations confirm these stability conditions. We monitor the quantity
\begin{equation}\label{eq:L2normTimeEv}
   \text{err-}a_{1} \coloneqq \dfrac{\left[ \int_{0}^{L_{y}} \left( a_{1}^{\mathrm{PML}}(x^{\ast},y,t) - a_{1}(x^{\ast},y,t)\right)^2 \dy \right]^{1/2}}{ \left[ \int_{0}^{L_{y}} \left( a_{1}(x^{\ast},y,t=0)\right)^2 \dy \right]^{1/2} },
\end{equation}
and plot its time evolution by performing a numerical simulation of the BGK+PML model with the same data described in Section \ref{sec:BGKPML_testing}.

Fig.~\ref{fig:err_a1_vs_time} reports on the results of the numerical simulations. In panel (a), we monitor the time evolution of err-$a_{1}$ for a set of $\LO$ values, namely $\left\lbrace 0, \, 10^{-4}, \, 10^{-3}, \, 10^{-2} \right\rbrace$. It appears that err-$a_{1}$ is minimal when $\LO$ is set to zero. 
From all the numerical simulations performed, it appears that we must set $\LO$ and $\LI$ to zero to address the stability issues and ensure an accurate description of the physical behavior. This agrees with the stability analysis performed in Section~\ref{subsec:mu4_CF}. We finally note that, even if we set $\LO$ to zero, the parameter $\AO$ still appears in the equations \eqref{eq:BGKPML_sigma2=0}.

Later, we will also analyze the sensitivity of our model to variations in the PML thickness $L$, so we proceed to establish some appropriate lower bound on $L$ before calculating the ANOVA expansion. Panel (b) of Fig.~\ref{fig:err_a1_vs_time} shows the time evolution of err-$a_{1}$, for several PML thicknesses, $L=\left\lbrace 0.10, \, 0.15, \, 0.20, \, 0.25 \right\rbrace$, with $\LO = 0$. A reasonable lower bound on the PML thickness $L$ is given by $L = 0.25$ because with this choice, err-$a_{1}$ is sufficiently small throughout the simulation.


\begin{figure}[htbp]
   \centering
   \subfloat[]{\includegraphics[width=0.475\textwidth]{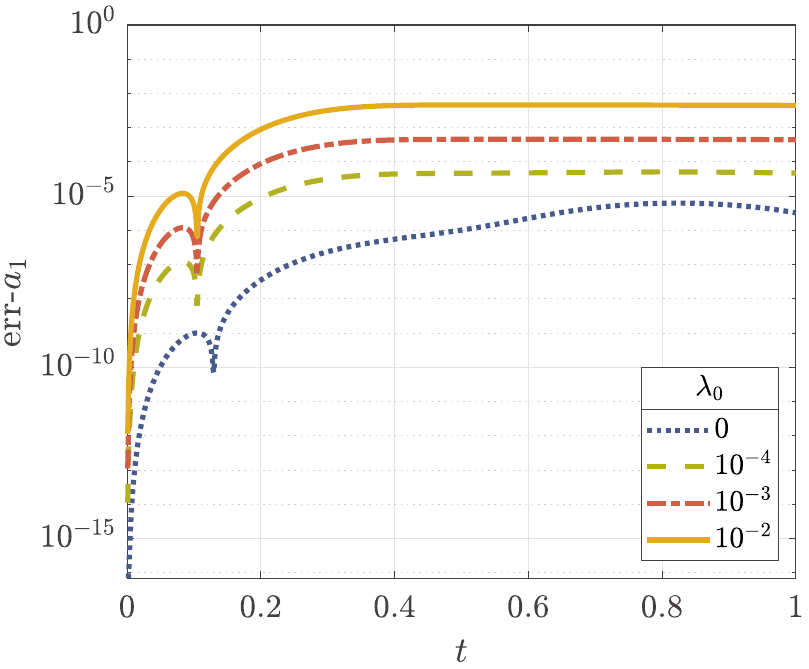}} \qquad
   \subfloat[]{\includegraphics[width=0.475\textwidth]{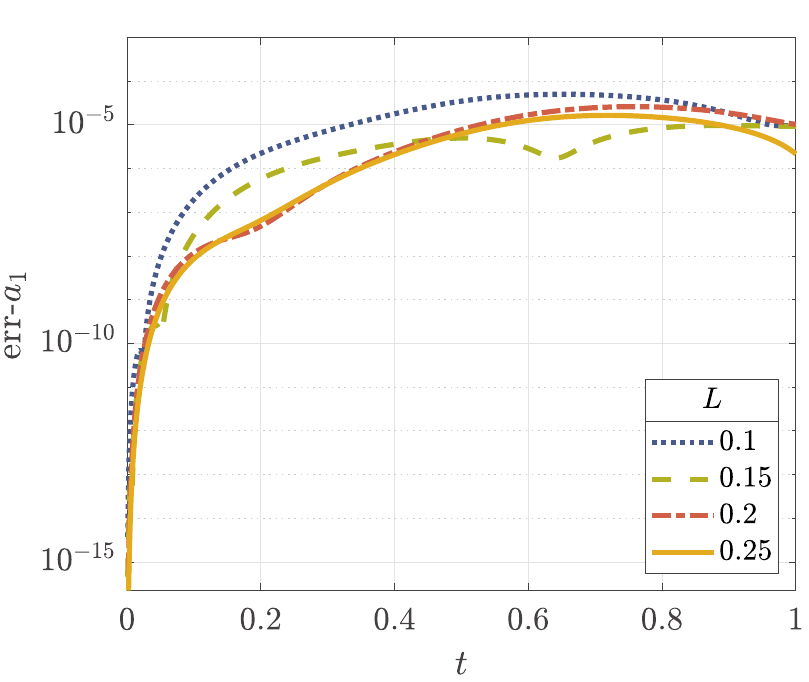}}
   \caption{Time evolution of err-$a_{1}$. Panel (a): For several values of $\LO$. Panel (b): For several PML thicknesses $L$.}   
   \label{fig:err_a1_vs_time}
\end{figure}

\subsection{ANOVA expansion applied to the BGK+PML model}\label{sec:BGKPML_ANOVA}

We now have all the ingredients to apply the ANOVA expansion to our functional~\eqref{eq:error_functional}.
We consider $g_1$ to be a function of $\AO$, $\AI$ and $L$, namely $g_1(\AO, \AI, L)$. We recall that we set $\LO$ and $\LI$ to zero to address the stability issues discussed in Sections~\ref{subsec:mu4_CF} and \ref{sec:LO_LI_to_zero}.

The TSIs are computed for $\beta=2,3,4$, being $\beta$ the exponent that appears in the expression~\eqref{eq:damping_function} of the damping function $\SI$. We let both $\AO$ and $\AI$ vary in the interval $[0,5]$, while the PML thickness $L$ take values in the interval $[0.25,0.80]$.
Table~\ref{tab:TSIs_AO_AI_L} reports the TSIs, obtained using increasingly accurate Gauss--Legendre cubature formulas, namely $(G_2)^3$, $(G_3)^3$ and $(G_4)^3$, for $\beta=2,3,4$.

\begin{table}[htbp]
   \caption{TSIs for the parameters $\AO$, $\AI$ and $L$,  using functional $g_1(\AO, \AI, L)$.}
   \label{tab:TSIs_AO_AI_L}
   \centering
   \begin{tabular}{ccccc}
                     PML exponent & Cubature type & {$\AO$} & {$\AI$} & {$L$}  \\ \toprule
   \multirow{3}{1.1cm}{$\beta=2$} & $(G_2)^3$     & 0.2274  & 0.2521  & 0.9435  \\
                                  & $(G_3)^3$     & 0.2251  & 0.2565  & 0.9478  \\
                                  & $(G_4)^3$     & 0.2221  & 0.2494  & 0.9607  \\ \midrule
   \multirow{3}{1.1cm}{$\beta=3$} & $(G_2)^3$     & 0.2212  & 0.2511  & 0.9586  \\
                                  & $(G_3)^3$     & 0.2104  & 0.2590  & 0.9716  \\
                                  & $(G_4)^3$     & 0.2112  & 0.2491  & 0.9705  \\ \midrule
   \multirow{3}{1.1cm}{$\beta=4$} & $(G_2)^3$     & 0.2201  & 0.2352  & 0.9640  \\
                                  & $(G_3)^3$     & 0.2114  & 0.2488  & 0.9814  \\
                                  & $(G_4)^3$     & 0.2051  & 0.2427  & 0.9865  \\ 
   \end{tabular}
\end{table}

Not surprisingly, the parameter with the largest TSI is the PML thickness $L$, while $\AO$ and $\AI$ basically have comparable sensitivity measures. Moreover, we observe that the TSIs are virtually independent of the PML exponent $\beta$.
We can also let $\beta$ be a parameter of the functional, i.e., we consider $g_1(\AO,\AI,\beta,L)$. We compute the TSIs assuming the following intervals of variation for the parameters
\[
   \AO \in \left[ 0, 3.5 \right], \quad \AI \in \left[ 0, 3.5 \right], \quad \beta \in \left[ 0, 4 \right], \quad L \in \left[ 0.25, 0.80 \right].
\]
The results are reported in Table~\ref{tab:g1_TSIsAO_AI_beta_L}. 

\begin{table}[htbp]
   \caption{TSIs for the parameters $\AO$, $\AI$, $\beta$ and $L$,  using functional $g_1(\AO, \AI, \beta, L)$.}
   \label{tab:g1_TSIsAO_AI_beta_L}
   \centering
   \begin{tabular}{ccccc}
      Cubature type & {$\AO$} & {$\AI$}   & {$\beta$}  & {$L$}  \\ \toprule
      $(G_2)^4$     & 0.1638  & 0.2474  & 0.2775 & 0.9312 \\
      $(G_3)^4$     & 0.1635  & 0.1916  & 0.2879 & 0.9385 \\ 
   \end{tabular}
\end{table}

As expected, the parameters $\AO$ and $\AI$ have only a minor impact on $g_1$. In contrast, the parameters $\beta$ and $L$ appearing in the definition of the damping function $\SI$ clearly dominate the functional outcome. In particular, the PML thickness $L$ has the most significant influence on $g_1$. This motivates us to investigate what happens if we set the least important parameters $\AO=\AI=1$, and we compute the ANOVA expansion of the functional $g_1(\beta,L)$. Table~\ref{tab:g1_TSIs_beta_L} reports the TSI values for this case. 

\begin{table}[htbp]
   \caption{TSIs for the parameters $\beta$ and $L$, using functional $g_1(\beta, L)$.}
   \label{tab:g1_TSIs_beta_L}
   \centering
   \begin{tabular}{ccc}
      Cubature type & {$\beta$}  & {$L$}    \\ \toprule
      $(G_2)^2$     & 0.0465     & 0.9535   \\
      $(G_3)^2$     & 0.0557     & 0.9443   \\
      $(G_4)^2$     & 0.0660     & 0.9340   \\ 
   \end{tabular}
\end{table}

The TSI values in Table~\ref{tab:g1_TSIs_beta_L} strongly confirm that the PML thickness $L$ is the parameter having the most significant influence on the PML behavior.

\subsection{Other functionals}\label{sec:other_functionals}

In this section, we consider other forms of the error functional to investigate whether and how the choice of the functional may influence the TSI values.

We define the error functional $g_2$ as
\[
g_2 = \dfrac{\int_{0}^{T} \left[ \int_{0}^{L_{y}} \left( a_{1}^{\mathrm{PML}}(x^{\ast},y,t) - a_{1}(x^{\ast},y,t)\right)^2 \dy \right]^{1/2} \mathrm{d}t}{ \left[ \int_{0}^{L_{y}} \left( a_{1}(x^{\ast},y,t=0)\right)^2 \dy \right]^{1/2} },
\]
Here, we are still looking at the error in the density $a_{1}$ on a vertical line close to the PML, but instead of taking the maximum of the $L^2$-norm, we compute its integral over the whole time interval; compare with \eqref{eq:error_functional}.
Table~\ref{tab:g2_TSIsAO_AI_beta_L} reports the values of the TSIs obtained when using the functional $g_2(\AO, \AI, \beta, L)$.

\begin{table}[htbp]
   \caption{TSIs for the parameters $\AO$, $\AI$, $\beta$ and $L$, using functional $g_2(\AO, \AI, \beta, L)$.}
   \label{tab:g2_TSIsAO_AI_beta_L}
   \centering
   \begin{tabular}{ccccc}
      Cubature type & {$\AO$} & {$\AI$}   & {$\beta$}  & {$L$}  \\ \toprule
      $(G_2)^4$     & 0.1623  & 0.1875    & 0.3061     & 0.9343 \\
      $(G_3)^4$     & 0.1596  & 0.1659    & 0.3740     & 0.9298 \\ 
   \end{tabular}
\end{table}

It is encouraging to observe that these values are very similar to those in Table~\ref{tab:g1_TSIsAO_AI_beta_L}, meaning that the model's sensitivity to the various parameters is virtually independent of the choice of the functional. The only noticeable effect of choosing $g_1$ over $g_2$ is that the TSI of $\AI$ decreases a bit while the TSI of $\beta$ increases.

Freezing $\AO=\AI=1$ and recomputing the ANOVA analysis using the functional $g_2(\beta, L)$, we obtain the results in Table~\ref{tab:g2_TSIs_beta_L}.

\begin{table}[htbp]
   \caption{TSIs for the parameters $\beta$ and $L$, using functional $g_2(\beta, L)$.}
   \label{tab:g2_TSIs_beta_L}
   \centering
   \begin{tabular}{ccc}
      Cubature type & {$\beta$}  & {$L$}    \\ \toprule
      $(G_2)^2$     & 0.0543     & 0.9457   \\
      $(G_3)^2$     & 0.0803     & 0.9197   \\
      $(G_4)^2$     & 0.0924     & 0.9076   \\ 
   \end{tabular}
\end{table}

Another option for the error functional is to calculate the $L^2$-norm of the error in $a_{1}$ on the entire domain $\Omega$ and then integrate over time, normalized with respect to the $L^{2}$-norm of the initial condition of $a_{1}$ over $\Omega$, i.e.,
\[
   g_3 = \dfrac{\int_{0}^{T} \left[ \int_{0}^{L_{x}}\int_{0}^{L_{y}} \left( a_{1}^{\mathrm{PML}}(x,y,t) - a_{1}(x,y,t)\right)^2 \dx \dy \right]^{1/2} \mathrm{d}t}{ \left[ \int_{0}^{L_{x}}\int_{0}^{L_{y}} \left( a_{1}(x,y,t=0)\right)^2 \dx \dy \right]^{1/2} }.
\]
This choice of functional yields the TSI values reported in Table~\ref{tab:g3_TSIsAO_AI_beta_L}, which again confirms that the TSIs are basically independent of the choice of error functional.
The results in Table~\ref{tab:g3_TSIs_beta_L}, obtained for the functional $g_3(\beta, L)$, also lead to similar observations.

\begin{table}[htbp]
   \caption{TSIs for the parameters $\AO$, $\AI$, $\beta$ and $L$, using functional $g_3(\AO, \AI, \beta, L)$.}
   \label{tab:g3_TSIsAO_AI_beta_L}
   \centering
   \begin{tabular}{ccccc}
      Cubature type & {$\AO$} & {$\AI$}   & {$\beta$}  & {$L$}   \\ \toprule
      $(G_2)^4$     & 0.1649  & 0.1773    & 0.3072     & 0.9374  \\
      $(G_3)^4$     & 0.1605  & 0.1627    & 0.3920     & 0.9280  \\ 
   \end{tabular}
\end{table}

\begin{table}[htbp]
   \caption{TSIs for the parameters $\beta$ and $L$, using functional $g_3(\beta, L)$.}
   \label{tab:g3_TSIs_beta_L}
   \centering
   \begin{tabular}{ccc}
      Cubature type & {$\beta$}  & {$L$}    \\ \toprule
      $(G_2)^2$     & 0.0533     & 0.9467   \\
      $(G_3)^2$     & 0.0845     & 0.9155   \\
      $(G_4)^2$     & 0.1053     & 0.8947   \\ 
   \end{tabular}
\end{table}

Summarizing, we have answered the question \emph{``What are the most important parameters in the BGK+PML model?''}. For stability issues, we know that we have to set $\LO$ and $\LI$ to zero, while the results of the sensitivity analysis show that the most significant parameters are the PML exponent $\beta$ and the PML thickness $L$. Moreover, we have seen that the values of the TSIs do not significantly depend on the choice of error functional.

\subsection{Isentropic vortex test case}

As an additional numerical experiment to investigate how the ANOVA results change, we analyze the isentropic vortex test case considered in~\protect{\cite{Hu:2008,Karakus:2019}}. Here, we only state those concepts necessary to understand the simulation; we refer the reader to~\protect{\cite{Hu:2008,Karakus:2019}} for more details. The two-dimensional nonlinear Euler equations support an advective solution of the form
\begin{equation}\label{eq:isentropic_vortex}
   \begin{pmatrix}
      \rho(\xbold, t)   \\
         u(\xbold, t)   \\
         v(\xbold, t)   \\
         p(\xbold, t)
   \end{pmatrix} =
   \begin{pmatrix}
         0   \\
         U_{0}   \\
         V_{0}   \\
         0
   \end{pmatrix} +
   \begin{pmatrix}
      \rho_{r}(r)   \\
      -u_{r}(r) \sin\theta   \\
      u_{r}(r) \cos\theta   \\
         p_{r}(r)
   \end{pmatrix},
\end{equation}
where $ (U_{0},V_{0}) $ is the constant advective velocity, $u_{r}(r) $ is the radial velocity distribution, and $ r = \sqrt{(x-U_{0} t)^{2} + (y-V_{0} t)^{2}} $. The solution \eqref{eq:isentropic_vortex} advects with constant velocity $ (U_{0},V_{0}) $. For our numerical experiments, we consider a radial velocity distribution of the form~\protect{\cite[(55)]{Karakus:2019}}
\begin{equation}\label{eq:radial_vel_vortex}
   u_{r}(r) = \frac{U_{\mathrm{max}}}{b} \, r \, e^{\frac{1}{2}\left(1-\frac{r^{2}}{b^{2}}\right)},
\end{equation}
where $ U_{\mathrm{max}} $ is the maximum velocity at $r = b$. The density distribution is~\protect{\cite[(55)]{Hu:2008}}
\begin{equation}\label{eq:density_distrib_vortex}
   \rho_{r}(r) = \left( 1 - \frac{1}{2} (\tau - 1) \, U_{\mathrm{max}}^{2} \, e^{1-\frac{r^{2}}{b^{2}}} \right)^{1/(\tau-1)}.
\end{equation}
The initial condition is the one given by~\eqref{eq:isentropic_vortex}, \eqref{eq:radial_vel_vortex}, and~\eqref{eq:density_distrib_vortex} with $t=0$, $ (U_{0},V_{0}) = (0.5, 0) $, $ U_{\mathrm{max}} = 0.5\,U_{0} = 0.25 $, and $ b = 0.2 $. The initial conditions for the BGK coefficients are then given by~\eqref{eq:a0a1a2} and~\eqref{eq:a3a4a5}. The physical domain of interest is the square $ [-1,1]^{2} $.
The BGK model without the PML is solved on a larger computational domain to let the vortex freely propagate. The absorption profile is taken as in~\eqref{eq:damping_function}, with width $ L = 0.50 $ and exponent $\beta = 4$. In contrast to~\protect{\cite[Sec. 4.1]{Hu:2008}} and~\protect{\cite[Sec. 5.1.1]{Karakus:2019}}, here we consider a PML only in the $x$-direction at the outflow boundary (i.e., the PML corresponds to the domain $ [1, 1.5] \times [-1,1] $). The computational domains of both the BGK and the BGK+PML models are discretized with a uniform grid with $ h_{x} = h_{y} = 0.01 $.

Fig.~\ref{fig:vv_contours_Vortex_BGKPML} shows the $v$-velocity contours at $t=0$, 1.5, 2.5, and 3.5, demonstrating the absorption of the vortex by the PML at the outflow boundary. Notably, the vortex also preserves its symmetry while penetrating the absorbing layer, a signal of minimal reflections at the interface. Overall, this figure provides qualitative evidence of the BGK+PML model effectiveness for the isentropic vortex test case, comparable to~\protect{\cite[Fig. 4]{Hu:2008}} and~\protect{\cite[Fig. 3]{Karakus:2019}}.


\begin{figure}[htbp]
   \centering
   \includegraphics[width=0.75\textwidth]{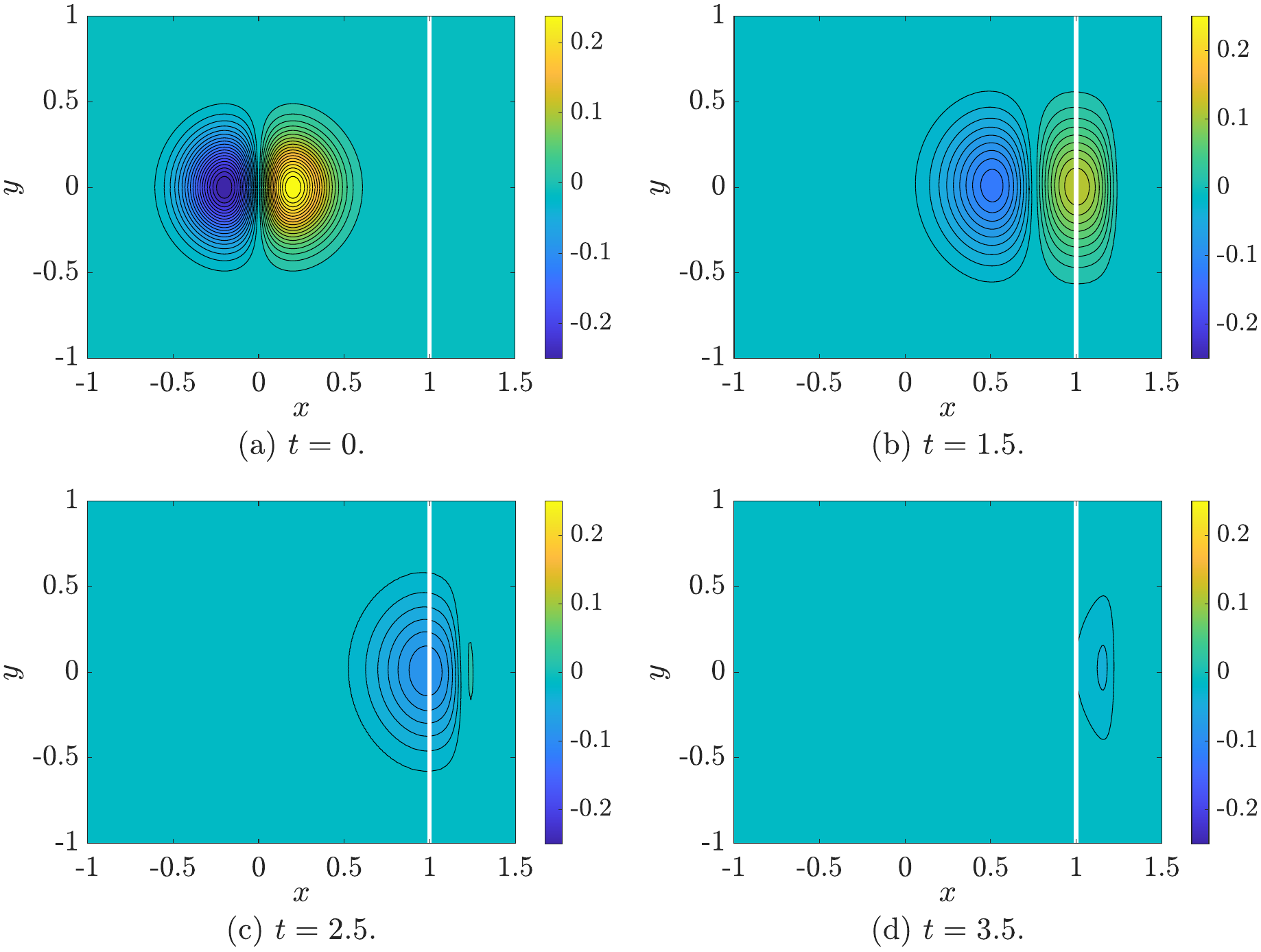}
   \caption{BGK+PML model for the isentropic vortex test case. Contours of the $v$-velocity distribution from $-0.25$ to $0.25$ with an increment of $0.0125$ excluding the zero level, at four different time instants.}
   \label{fig:vv_contours_Vortex_BGKPML}
\end{figure}

For computing the TSIs of the BGK+PML model for the current test case, we use a coarser grid, with a uniform discretization parameter $ h_{x} = h_{y} = 0.1 $, so that a larger time step $ \Delta t = 0.025 $ can be used for the numerical time integration.

As we did for the first validation test case in Sections~\ref{sec:BGKPML_ANOVA} and~\ref{sec:other_functionals}, we consider the difference between the numerical solution of the BGK+PML model and a reference solution obtained using the plain BGK model on a larger computational domain, along a vertical line close to the outflow boundary, $ x^{\ast} = 0.9 $.
We first consider the functional
\[
   h_{1} = \displaystyle\max_{t\in \left[ 0,T \right]} \left\lbrace \left\lVert v^{\mathrm{PML}}(x^{\ast},y,t) - v(x^{\ast},y,t) \right\rVert_{\infty} \right\rbrace,
\]
where $T=3.5$, and we assume the following intervals of variation for the parameters
\[
   \AO \in \left[ 0, 3.5 \right], \quad \AI \in \left[ 0, 3.5 \right], \quad \beta \in \left[ 0, 4 \right], \quad L \in \left[ 0.1, 1 \right].
\]
Clearly, $h_{1}$ is calculated on the physical domain only, namely for $x^{\ast} = 0.9$ and $ y \in [-1,1] $. We perform the ANOVA expansion and compute the TSIs as explained in the previous sections. Table~\ref{tab:h1_AO_AI_beta_L_vortex} reports on the numerical results.

\begin{table}[htbp]
   \caption{TSIs for the parameters $\AO$, $\AI$, $\beta$ and $L$, using functional $g_1(\AO, \AI, \beta, L)$, for the isentropic vortex test case.}
   \label{tab:h1_AO_AI_beta_L_vortex}
   \centering
   \begin{tabular}{*{5}{c}}
      Cubature type   &  {$\AO$}  &  {$\AI$}  &  {$\beta$}  &    {$L$}  \\ \toprule
      $(G_{2})^{4}$   &  0.2026   &   0.2323  &    0.3066   &   0.7753  \\
      $(G_{3})^{4}$   &  0.2153   &   0.2783  &    0.3280   &   0.7431
   \end{tabular}
\end{table}
            
As in Sections~\ref{sec:BGKPML_ANOVA} and~\ref{sec:other_functionals}, we investigate what happens when freezing $\AO=\AI=1$ and redoing the ANOVA expansion using the functional $h_{1}(\beta, L)$, we obtain the TSIs reported in Table~\ref{tab:h1_beta_L_vortex}, for increasing precision of the numerical integration formulas.

\begin{table}[htbp]
   \caption{TSIs for the parameters $\beta$ and $L$, using functional $h_{1}(\beta, L)$, for the isentropic vortex test case.}
   \label{tab:h1_beta_L_vortex}
   \centering
   \begin{tabular}{*{3}{c}}
      Cubature type   &  {$\beta$}  &    {$L$}  \\ \toprule
      $(G_{2})^{2}$   &    0.1480   &   0.8520  \\
      $(G_{3})^{2}$   &    0.1503   &   0.8497  \\
      $(G_{4})^{2}$   &    0.1547   &   0.8453
   \end{tabular}
\end{table}

As in the previous test case, the impact of $\alpha_{0}$ and $\alpha_{1}$ is limited. The parameters $\beta$ and $L$ defining the damping function $\sigma_{1}$ have the highest TSI values.

Another option for the error functional is to calculate the $L^2$-norm of the error in the $v$-velocity on the physical domain $ [-1,1]^{2} $ and then integrate over time, normalized with respect to the $L^{2}$-norm of the initial condition of $v$ over $ [-1,1]^{2} $, i.e.,
\[
   h_{2} = \dfrac{\int_{0}^{T} \left[ \int_{-1}^{1}\int_{-1}^{1} \left( v^{\mathrm{PML}}(x,y,t) - v(x,y,t)\right)^2 \dx \dy \right]^{1/2} \mathrm{d}t}{ \left[ \int_{-1}^{1}\int_{-1}^{1} \left( v(x,y,t=0)\right)^2 \dx \dy \right]^{1/2} },
\]
where $T = 3.5$. Table~\ref{tab:h2_AO_AI_beta_L_vortex} reports the results for the TSIs when using $ h_{2} $ as a function of the parameters $\AO$, $\AI$, $\beta$ and $L$.

\begin{table}[htbp]
   \caption{TSIs for the parameters $\AO$, $\AI$, $\beta$ and $L$, using functional $h_{2}(\AO, \AI, \beta, L)$, for the isentropic vortex test case.}
   \label{tab:h2_AO_AI_beta_L_vortex}
   \centering
   \begin{tabular}{*{5}{c}}
      Cubature type   &  {$\AO$}  &  {$\AI$}  &  {$\beta$}  &    {$L$}  \\ \toprule
      $(G_{2})^{4}$   &  0.1911   &   0.1939  &    0.3045   &   0.8051  \\
      $(G_{3})^{4}$   &  0.1914   &   0.1997  &    0.3448   &   0.7831
   \end{tabular}
\end{table}
            
Freezing $\AO=\AI=1$ and redoing the ANOVA expansion using the functional $h_{2}(\beta, L)$, we obtain the TSIs reported in Table~\ref{tab:h2_beta_L_vortex}. As expected, as we increase the accuracy of the Gauss--Legendre integration formula, the TSI values in Tables~\ref{tab:h2_AO_AI_beta_L_vortex} and~\ref{tab:h2_beta_L_vortex} seem to converge.

\begin{table}[htbp]
   \caption{TSIs for the parameters $\beta$ and $L$, using functional $h_{2}(\beta, L)$, for the isentropic vortex test case.}
   \label{tab:h2_beta_L_vortex}
   \centering
   \begin{tabular}{*{3}{c}}
      Cubature type   &  {$\beta$}  &    {$L$}   \\ \toprule
      $(G_{2})^{2}$   &    0.1548   &   0.8452   \\
      $(G_{3})^{2}$   &    0.1646   &   0.8354   \\
      $(G_{4})^{2}$   &    0.1651   &   0.8349
   \end{tabular}
\end{table}

Comparing Table~\ref{tab:h2_AO_AI_beta_L_vortex} to Table~\ref{tab:h1_AO_AI_beta_L_vortex}, and Table~\ref{tab:h2_beta_L_vortex} to Table~\ref{tab:h1_beta_L_vortex} demonstrates that the TSIs are basically independent of the choice of the functional, which confirms the findings for the previous test case. Also, the TSI values themselves are not too far from those in Sections~\ref{sec:BGKPML_ANOVA} and~\ref{sec:other_functionals}. The main difference we notice is that for the isentropic vortex test case, the TSIs related to $\beta$ are slightly higher. The range of TSI values for $ L $ in the first test case is the interval $ [0.89, 0.99] $, while for the isentropic vortex is $ [0.74, 0.85] $. Indeed, if we focus our attention on those functionals that only depend on $\beta$ and $L$ and look at their TSIs, we observe that for the first test case, the TSI for $ \beta $ take values in the range $ [ 0.0465, 0.1053 ] $, while for the second test case, the interval of variation is $ [ 0.1480, 0.1651 ] $. Analogously, the TSIs for $L$ take values in $ [ 0.8947, 0.9535 ] $ for the first test case and in $ [ 0.8349, 0.8520 ] $ for the second test case.

\section{Conclusions and outlook}\label{sec:conclusions}

In this work, we studied the stability and sensitivity of an absorbing PML layer for the BGK approximation of the Boltzmann equation. In particular, we investigated the parameters' role and significance in the BGK+PML model of \protect{\cite{Gao2011a}}. To establish reasonable parameter bounds, we analyzed the model's stability, focusing on the symbol of the differential operator of the system, utilizing the energy decay and a technique based on continued fractions. Our analyses reveal that setting $\LO$ and $\LI$ to zero is necessary to ensure stability, which is also confirmed by extensive numerical simulations. Additionally, our analyses provided some bounds for the other parameters in the BGK+PML model.

We used the ANOVA expansion of multivariate functions to calculate the total sensitivity indices (TSIs) of the parameters for two different test cases and for several choices of error functional. The parameter values required for the simulations to obtain the TSIs were carefully chosen based on the bounds established through the stability analysis. Via the sensitivity analysis, we identified the most crucial parameters, namely, the PML exponent $\beta$ and thickness $L$. Additionally, the numerical experiments demonstrate that the TSIs are essentially independent of the choice of error functional.

As a future research outlook, we emphasize that our study might be biased toward a limited test set and implementation details. One could change the problem configuration, particularly the initial and boundary conditions, and explore these aspects more systematically. For instance, in this study, we always considered an initial value problem, but we could also consider a continuously excited problem. The discretization method could also be changed, e.g., one could use the discontinuous Galerkin method instead of finite differences.
We finally recall that, for low Mach numbers and weakly compressible flows, one could couple the Navier--Stokes equations (NSE) and the BGK model by solving the NSE in the physical domain and the BGK model in the PML domain.
The coupling of the NSE and the BGK equations is an open research direction which is left for future investigation.

\section*{Acknowledgments}

The first author started this work as a master's thesis at the EPFL under the supervision of Jan S. Hesthaven. It was completed during the first author's postdoctoral fellowship at the National Center for Theoretical Sciences in Taiwan (R.O.C.).

\section*{Data availability}

The code and datasets generated and analyzed during the current study are available upon reasonable request.

\section*{Conflict of interest}

The authors declare that they have no conflict of interest.

\bibliographystyle{aomalpha}

\begin{footnotesize}
   \bibliography{PML_preprint_arXiv_2024_biblio.bib}

\providecommand{\bysame}{\leavevmode\hbox to3em{\hrulefill}\thinspace}
\providecommand{\noopsort}[1]{}
\providecommand{\mr}[1]{\href{http://www.ams.org/mathscinet-getitem?mr=#1}{MR~#1}}
\providecommand{\zbl}[1]{\href{http://www.zentralblatt-math.org/zmath/en/search/?q=an:#1}{Zbl~#1}}
\providecommand{\jfm}[1]{\href{http://www.emis.de/cgi-bin/JFM-item?#1}{JFM~#1}}
\providecommand{\arxiv}[1]{\href{http://www.arxiv.org/abs/#1}{arXiv~#1}}
\providecommand{\doi}[1]{\url{https://doi.org/#1}}
\providecommand{\MR}{\relax\ifhmode\unskip\space\fi MR }
\providecommand{\MRhref}[2]{%
  \href{http://www.ams.org/mathscinet-getitem?mr=#1}{#2}
}
\providecommand{\href}[2]{#2}
\begin{thebibliography}{KCHW19}

\bibitem[And97]{Andres1997}
\bgroup\scshape{}T.~Andres\egroup{}, {Sampling methods and sensitivity analysis
  for large parameter sets},  \emph{J. Stat. Comput. Sim.} \textbf{57} no.~1-4
  (1997), 77--110. \doi{https://doi.org/10.1080/00949659708811804}.
\bibitem[AHK06]{Appelo2006}
\bgroup\scshape{}D.~Appelö\egroup{}, \bgroup\scshape{}T.~Hagstrom\egroup{},
  and \bgroup\scshape{}G.~Kreiss\egroup{}, {Perfectly matched layers for
  hyperbolic systems: general formulation, well-posedness, and stability},
  \emph{SIAM J. Appl. Math.} \textbf{67} no.~1 (2006), 1--23.
  \doi{https://doi.org/10.1137/050639107}.
\bibitem[Ber94]{Berenger1994}
\bgroup\scshape{}J.-P. Berenger\egroup{}, {A perfectly matched layer for the
  absorption of electromagnetic waves},  \emph{J. Comput. Phys.} \textbf{114}
  no.~2 (1994), 185--200. \doi{https://doi.org/10.1006/jcph.1994.1159}.
\bibitem[BGK54]{BGK:1954}
\bgroup\scshape{}P.~L. Bhatnagar\egroup{}, \bgroup\scshape{}E.~P.
  Gross\egroup{}, and \bgroup\scshape{}M.~Krook\egroup{}, {A Model for
  Collision Processes in Gases. I. Small Amplitude Processes in Charged and
  Neutral One-Component Systems},  \emph{Phys. Rev.} \textbf{94} no.~3 (1954),
  511--525. \doi{https://link.aps.org/doi/10.1103/PhysRev.94.511}.
\bibitem[CCG09]{Cao2009}
\bgroup\scshape{}Y.~Cao\egroup{}, \bgroup\scshape{}Z.~Chen\egroup{}, and
  \bgroup\scshape{}M.~Gunzburger\egroup{}, {ANOVA expansions and efficient
  sampling methods for parameter dependent nonlinear PDEs},  \emph{Int. J.
  Numer. Anal. Model.} \textbf{6} no.~2 (2009), 256--273.
\bibitem[CC95]{Chapman_Cowling:1995}
\bgroup\scshape{}S.~Chapman\egroup{} and \bgroup\scshape{}T.~G.
  Cowling\egroup{}, \emph{The mathematical theory of non-uniform gases}, third
  ed., \emph{Cambridge Math. Lib.}, Cambridge University Press, Cambridge,
  1995, An account of the kinetic theory of viscosity, thermal conduction and
  diffusion in gases, In co-operation with D. Burnett, With a foreword by Carlo
  Cercignani. \mr{1148892}.
\bibitem[Eva10]{evans2010partial}
\bgroup\scshape{}L.~C. Evans\egroup{}, \emph{{Partial Differential Equations}},
  second ed., \emph{Grad. Stud. Math.} \textbf{19}, American Mathematical
  Society, Providence, RI, March 2010. \doi{https://doi.org/10.1090/gsm/019}.
\bibitem[Fra46]{Frank1946}
\bgroup\scshape{}E.~Frank\egroup{}, {On the zeros of polynomials with complex
  coefficients},  \emph{Bull. Amer. Math. Soc.} \textbf{52} no.~2 (1946),
  144--157. \doi{https://doi.org/10.1090/s0002-9904-1946-08526-2}.
\bibitem[GH11]{Gao2011}
\bgroup\scshape{}Z.~Gao\egroup{} and \bgroup\scshape{}J.~S. Hesthaven\egroup{},
  {Efficient solution of ordinary differential equations with high-dimensional
  parametrized uncertainty},  \emph{Commun. Comput. Phys.} \textbf{10} no.~2
  (2011), 253--278. \doi{https://doi.org/10.4208/cicp.090110.080910a}.
\bibitem[GHW11]{Gao2011a}
\bgroup\scshape{}Z.~Gao\egroup{}, \bgroup\scshape{}J.~S. Hesthaven\egroup{},
  and \bgroup\scshape{}T.~Warburton\egroup{}, \emph{{Efficient absorbing layers
  for weakly compressible flows}}, Submitted to Journal of scientific computing
  (Springer, ISSN 0885-7474), 2011.
\bibitem[GST01]{Gottlieb2001}
\bgroup\scshape{}S.~Gottlieb\egroup{}, \bgroup\scshape{}C.-W. Shu\egroup{}, and
  \bgroup\scshape{}E.~Tadmor\egroup{}, {Strong stability-preserving high-order
  time discretization methods},  \emph{SIAM Rev.} \textbf{43} no.~1 (2001),
  89--112. \doi{https://doi.org/10.1137/s003614450036757x}.
\bibitem[Gra49]{Grad1949}
\bgroup\scshape{}H.~Grad\egroup{}, {On the kinetic theory of rarefied gases},
  \emph{Comm. Pure Appl. Math.} \textbf{2} no.~4 (1949), 331--407.
  \doi{https://doi.org/10.1002/cpa.3160020403}.
\bibitem[GKO13]{Gustafsson1995}
\bgroup\scshape{}B.~Gustafsson\egroup{}, \bgroup\scshape{}H.-O.
  Kreiss\egroup{}, and \bgroup\scshape{}J.~Oliger\egroup{},
  \emph{{Time-Dependent Problems and Difference Methods}}, second ed.,
  \emph{Pure and Appl. Math. (Hoboken)}, John Wiley \& Sons, Inc., Hoboken, NJ,
  September 2013. \doi{10.1002/9781118548448}.
\bibitem[Hag03]{Hagstrom2003}
\bgroup\scshape{}T.~Hagstrom\egroup{}, {A new construction of perfectly matched
  layers for hyperbolic systems with applications to the linearized Euler
  equations},  \emph{Mathematical and numerical aspects of wave
  propagation---{WAVES} 2003} (2003), 125--129.
  \doi{https://doi.org/10.1007/978-3-642-55856-6_20}.
\bibitem[Her87]{Hernquist:1987}
\bgroup\scshape{}D.~C. Hernquist\egroup{}, {Smoothly symmetrizable hyperbolic
  systems of partial differential equations},  \emph{Math. Scand.} \textbf{61}
  no.~2 (1987), 262--275. \doi{https://doi.org/10.7146/math.scand.a-12203}.
\bibitem[HLL08]{Hu:2008}
\bgroup\scshape{}F.~Q. Hu\egroup{}, \bgroup\scshape{}X.~Li\egroup{}, and
  \bgroup\scshape{}D.~Lin\egroup{}, {Absorbing boundary conditions for
  nonlinear Euler and Navier--Stokes equations based on the perfectly matched
  layer technique},  \emph{J. Comput. Phys.} \textbf{227} no.~9 (2008),
  4398--4424. \doi{https://doi.org/10.1016/j.jcp.2008.01.010}.
\bibitem[KCHW19]{Karakus:2019}
\bgroup\scshape{}A.~Karakus\egroup{}, \bgroup\scshape{}N.~Chalmers\egroup{},
  \bgroup\scshape{}J.~S. Hesthaven\egroup{}, and
  \bgroup\scshape{}T.~Warburton\egroup{}, {Discontinuous Galerkin
  discretizations of the Boltzmann--BGK equations for nearly incompressible
  flows: semi-analytic time stepping and absorbing boundary layers},  \emph{J.
  Comput. Phys.} \textbf{390} (2019), 175--202.
  \doi{https://doi.org/10.1016/j.jcp.2019.03.050}.
\bibitem[LeV07]{LeVeque2007}
\bgroup\scshape{}R.~J. LeVeque\egroup{}, \emph{{Finite Difference Methods for
  Ordinary and Partial Differential Equations}}, Society for Industrial and
  Applied Mathematics (SIAM), Philadelphia, PA, January 2007.
  \doi{https://doi.org/10.1137/1.9780898717839}.
\bibitem[Mar66]{Marden1966}
\bgroup\scshape{}M.~Marden\egroup{}, \emph{{Geometry of Polynomials}},
  \emph{Math. Surveys} \textbf{No. 3}, American Mathematical Society,
  Providence, RI, December 1966. \doi{https://doi.org/10.1090/surv/003}.
\bibitem[Rez05]{Rezzolla2005}
\bgroup\scshape{}L.~Rezzolla\egroup{}, {Numerical Methods for the Solution of
  Hyperbolic Partial Differential Equations},  (2005).
\bibitem[SCS00]{Saltelli2000}
\bgroup\scshape{}A.~Saltelli\egroup{}, \bgroup\scshape{}K.~Chan\egroup{}, and
  \bgroup\scshape{}E.~M. Scott\egroup{} (eds.), \emph{{Sensitivity Analysis}},
  \emph{Wiley Ser. Probab. Stat.}, Wiley, Chichester, 2000.
\bibitem[Sut15]{Sutti:2015}
\bgroup\scshape{}M.~Sutti\egroup{}, \emph{{Analysis and Optimization of
  Perfectly Matched Layers for the Boltzmann Equation}}, Master's thesis, EPFL,
  Computational Mathematics and Simulation Science, June 2015.
  \doi{http://refhub.elsevier.com/S0021-9991(19)30233-5/bib73757474695F616E616C797369735F32303135s1}.
\bibitem[TKSR00]{Tolke:2000}
\bgroup\scshape{}J.~Tölke\egroup{}, \bgroup\scshape{}M.~Krafczyk\egroup{},
  \bgroup\scshape{}M.~Schulz\egroup{}, and \bgroup\scshape{}E.~Rank\egroup{},
  {Discretization of the Boltzmann equation in velocity space using a Galerkin
  approach},  \emph{Comput. Phys. Comm.} \textbf{129} no.~1 (2000), 91--99.
  \doi{https://doi.org/10.1016/S0010-4655(00)00096-5}.
\bibitem[WF03]{Wang:2003}
\bgroup\scshape{}X.~Wang\egroup{} and \bgroup\scshape{}K.-T. Fang\egroup{},
  {The effective dimension and quasi-Monte Carlo integration},  \emph{J.
  Complex.} \textbf{19} no.~2 (2003), 101--124.
  \doi{https://doi.org/10.1016/S0885-064X(03)00003-7}.
\end{thebibliography}
\end{footnotesize}

\end{document}